\documentclass{amsart}[12pt]
\usepackage{a4wide}
\usepackage[polish]{babel}
\begin{document}
\sloppy
\title{Geometria nieortodoksyjna}
\author{Piotr "Sniady}
\thanks{Pierwsza wersja tej
ksi"a"reczki powsta"la jako przewodnik po
warsztatach o krzywiznie, kt"ore autor prowadzi"l na obozie stypendyst"ow
Krajowego Funduszu na rzecz Dzieci w Jadwisinie 12--24 kwietnia 1997.
Adres do korespondencji: ul. "Sliczna 49/12 50-566~Wroc"law. 
E-mail:{\tt psnia@math.uni.wroc.pl}}

\maketitle

\renewcommand{\tan}{\rm{tg\ }}
\newcommand{\rozwiazanie}[1]{\vspace{1ex}{\Large \bf #1}}
\newenvironment{lista}{\begin{list}{$\bullet$}{\itemsep0ex \topsep0ex
\parsep0ex}}{\end{list}}
\newcommand{\wektor}[1]{\stackrel{\longrightarrow}{#1}}

\section{Wst"ep}
Wsta"n i zr"ob krok przed siebie. Obr"o"c si"e w prawo o $90$ stopni
i zn"ow zr"ob krok przed siebie. Jeszcze raz. I jeszcze raz. Gdzie jeste"s 
teraz?
To proste. Chodzi"le"s po kwadracie i teraz zn"ow jeste"s tam, gdzie na
pocz"atku. Ale czy na pewno?
Co si"e stanie, je"sli b"edziesz robi"c dwa razy takie kroki?
Wydaje si"e, "re b"edziesz chodzi"c po dwa razy wi"ekszym kwadracie
i teraz te"r wr"ocisz
do punktu pocz"atkowego. A co je"sli b"edziesz robi"c naprawd"e
{\Huge du\-"re kroki}? Naprawd"e du"re, na przyk"lad d"lugo"sci $10000$
kilometr"ow? Zobaczmy...

Wcale nie wr"ocili"smy do punktu wyj"scia! Czy nie jest to zastanawiaj"ace?
Co ciekawsze, do punktu wyj"scia wr"ocili"smy ju"r po wykonaniu trzech krok"ow.
Nasze stopy wykre"sli"ly boki tr"ojk"ata, kt"orego {\bf wszystkie k"aty
s"a proste}! Zatem wbrew temu, czego uczy si"e w szko"lach,
{\bf suma  k"at"ow w naszym tr"ojk"acie nie jest
r"owna 180 stopni}!

No jasne, zadaj"ac pytanie o chodzeniu i zakr"ecaniu nie doda"lem, "re
mamy chodzi"c po powierzchni kuli ziemskiej. Gdyby"smy chodzili po p"laszczy"znie,
takich dziwnych zjawisk by"smy na pewno nie znale"zli. Jaka cecha r"o"rni
wi"ec p"laszczyzn"e i sfer"e? Jest ni"a {\bf krzywizna}. I
w"la"snie
badaniu wspania"lych perspektyw, jakie ona nam daje, jest po"swi"econa ta ksi"a"reczka.

Na pocz"atku zbadamy jak wygl"ada uprawianie geometrii na najprostszej
powierzchni obdarzonej krzywizn"a--na sferze.
Spr"obujemy udowodni"c twierdzenie Pitagorasa i nauczymy si"e
do"s"c dziwnych sposob"ow na mierzenie pola powierzchni. No i u"sci"slimy, co to
jest krzywizna.
Potem zabierzemy si"e za do"s"c dowolne powierzchnie zakrzywione. Zastanowimy
si"e czym zast"api"c poj"ecie odcinka i prostej. B"edziemy starali si"e
g"l"ebiej zrozumie"c dlaczego id"ac prosto, idziemy tak krzywo.
Dowiemy si"e, co
to jest charakterystyka Eulera i udowodnimy niezwyk"ly wz"or na ca"lkowit"a
krzywizn"e np. fili"ranki.
Przy okazji mo"re zrozumiemy troch"e og"olnej teorii wzgl"edno"sci i
dlaczego wszystko spada na d"o"l.

Poniewa"r nie jest to podr"ecznik, ale raczej przeno"sny inspirator,
przygotowa"lem nieco problem"ow i
zada"n. S"a bardzo r"o"rne: niekt"ore s"a do"s"c proste, inne
naprawd"e trudne. Jest ich na  tyle du"ro, "re chyba nie da si"e ich
wszystkich rozwi"aza"c, zatem r"ownie"r od was zale"ry, kt"ore wybierzecie.
Rozdzia"l o spacerach na parkietach zawiera naprawd"e trudne problemy.
Je"sli chcesz napisa"c prac"e na Konkurs Prac Uczniowskich z Matematyki
i szukasz tematu...

Je\"sli spodoba wam si"e taka geometria, mo"re zechcecie poczyta"c co"s
wi"ecej? W razie jakich"s k"lopot"ow, mo"recie do mnie pisa"c.

Chcia"lbym przeprosi"c za brak jakichkolwiek rysunk"ow. {\em Jaki jest po"rytek z
ksi"a"rki, w kt"orej nie ma ani konwersacji, ani ilustracji?\footnote {Lewis Carroll, Alicja w krainie czar"ow.}}

{\bf Do niniejszych warsztat"ow napisa"lem modu"ly do pisania w"lasnych
prostych program"ow w Pascalu, kt"ore pozwol"a wam namacalnie
poeksperymentowa"c z krzywizn"a. Aby je otrzyma"c bez\-p"la\-tnie,
nale"ry przys"la"c do mnie pust"a dyskietk"e 3.5". Jedynym
warunkiem, jaki musisz spe"lni"c, by w pe"lni legalnie korzysta"c z
tych program"ow jest przes"lanie na m"oj adres "ladnej poczt"owki.
Zerknij do pliku {\tt czytaj.to!} przed u"ryciem.}

\section{Co warto czyta"c?}
Co {\bf na pewno} warto poczyta"c?
\begin{lista}
\item H. Abelson, A. A. diSessa, Geometria "r"o"lwia. {\em Bardzo
dobra ksi"a"rka, wbrew pozorom nie informatyczna, ale w"la"snie
matematyczna. Si"egnijcie po ni"a koniecznie!}
\item W. Krysicki, H. Pisarewska, T. "Swi"atkowski, Z geometri"a za pan brat.
\item M. Kordos, Pole tr"ojk"ata sferycznego i wz"or Eulera. Delta 5/1996.
\item M. Kordos, O r"o"rnych geometriach.
\item Feynmana wyk"lady z fizyki, tom II cz"e"s"c 2, rozdzia"l 42.
{\it Ten rozdzia"l jest o przestrzeni zakrzywionej, ale warto czyta"c
te"r inne rozdzia"ly. Jest to przepi"ekna ksi"a"rka.}
\item Matematyka wsp"o"lczesna, dwana"scie esej"ow. red. L. A. Steen,
esej Rogera Penrose'a.
\item G. Gamow, Mister Tompkins w krainie czar"ow. {\em Wspania"la
ba"s"n.}
\item H. M. S. Coxeter, Geometria dawna i nowa. {\em Polecam ka"rdemu
geometrze, cho"c akurat ma"lo jest w niej mowy o tym, czym tu si"e
zajmujemy}
\item S. Kulczycki, Geometria nieeuklidesowa.
\item D. Hilbert, S. Cohn-Vossen, Geometria pogladowa.
\end{lista}

\noindent
Poni"rsze pozycje s"a ma"lo ciekawe, ale jest w nich troch"e
wzork"ow trygonometrii sferycznej.
\begin{lista}
\item Encyklopednia szkolna, Matematyka, s. 290-291.
\item red. I. Dziubi"nski, T. "Swi"atkowski, Poradnik matematyczny, PWN, rozdzia"l V.2
\item Poradnik in"ryniera - Matematyka, rozdzia"l 2.
\end{lista}

\noindent
A je"sli jeste"s ju"r koneserem...
\begin{lista}
\item B. Schutz, Wst"ep do og"olnej teorii wzgl"edno"sci.
{\em Wbrew temu, co m"owi tytu"l (a mo"re zgodnie z tytu"lem?), jest
to wspania"ly, prosty podr"ecznik geometrii r"o"rniczkowej.}
\item Gallot, Riemannian geometry. {\em Pi"ekne, ale po angielsku.}
\item L. Landau, E. Lifszic, Teoria pola, rozdzia"ly 10 i 11. {\em Na
pocz"atek lepiej poczyta"c Schutza.}
\item K. Maurin, Analiza tom II. {\em Uwaga! Ta ksi"a"rka nie jest
podr"ecznikiem. Wielu ludzi, kt"orzy tak pr"obowali j"a czyta"c:
grzecznie, strona po stronie, zw"atpi"lo w swoje si"ly i wpad"lo w
kompleksy. Ale chyba warto zobaczy"c np. jak Maurin na kilka sposob"ow
definiuje wektory.}
\item Prawie dowolna ksi"a"rka z {\it geometri"a
r"o"rniczkow"a} w tytule.
\end{lista}

\section{Krzywa na p"laszczyznie.}
Popatrz na r"o"rne krzywe na p"laszczyznie. Na przyk"lad na prost"a.
Chyba wszyscy si"e zgodzimy, "re jest ona prosta i nie jest wcale skrzywiona.

A okr"ag? Tak, okr"ag na pewno jest skrzywiony. Im wi"ekszy okr"ag,
tym bardziej staje si"e podobny do prostej i tym mniej jest skrzywiony, a im
mniejszy--tym bardziej.

Maj"ac ju"r teraz jakie"s intuicje czym krzywizna jest, {\bf zdefiniujemy
krzywizn"e okr"egu} po prostu {\bf jako odwrotno"s"c jego promienia}
(a jaka jest krzywizna prostej?).

A co z krzywizn"a bardziej skomplikowanych krzywych? Jedno wida"c na
pierwszy rzut oka: mianowicie krzywizna zazwyczaj nie b"edzie sta"l"a
wie\-lko\-"sci"a na ca"lej d"lugo"sci krzywej, ale b"edzie zmienia"c
si"e od punktu do punktu.

Tak wi"ec, aby zmierzy"c krzywizn"e jakiej"s krzywej w pewnym
punkcie, spr"obujmy jak najdok"ladniej przybli"ry"c j"a w tym punkcie
przy pomocy okr"egu. Krzywizn"e krzywej przyjmiemy za r"own"a
krzywiznie tego okr"egu. Ale co to w"la"sciwie znaczy: jak najdok"ladniej
przybli"ry"c krzyw"a przy pomocy okr"egu?

Ja proponuj"e tak"a metod"e: {\bf wybieramy sobie dwa punkty krzywej
bardzo bliskie punktowi, w kt"orym chcemy liczy"c krzywizn"e. W ka"rdym
z tych dw"och punkt"ow rysujemy prost"a normaln"a do krzywej (normalna--to prosta prostopad"la do stycznej).
Punkt przeci"ecia
normalnych powinien by"c ju"r bardzo blisko od "srodka poszukiwanego
okr"egu.}

{\small Jeszcze tylko drobny problem techniczny. Mianowicie powinni"smy tak
otrzymany promie"n krzywizny i krzywizn"e opatrywa"c znakiem plus lub
minus zale"rnie od tego, w kt"or"a stron"e wygi"eta jest nasza krzywa.}

\subsection{Krzywizna paraboli i elipsy.}
\label{parabola}
Jaka jest krzywizna paraboli $y=a x^2$ mierzona w punkcie $(0,0)$?
Jaka jest krzywizna elipsy o d"lugiej p"o"losi $a$ i odleg"lo"sci
ognisk $2b$ mierzona w wierzcho"lkach?

\subsection{Ca"lkowita krzywizna krzywej.}
\label{calkowita}
Rozpatrzmy krzyw"a zamkni"et"a, kt"ora sk"lada si"e z kilku
"luk"ow okr"eg"ow tak dobranych, "reby nasza krzywa by"la g"ladka,
bez 'kant"ow'. Na ka"rdym "luku okr"egu krzywizna jest sta"la.
Dla ka"rdego "luku mno"rymy jego d"lugo"s"c i jego krzywizn"e.
Sum"e tych liczb nazywamy ca"lkowit"a krzywizn"a krzywej.

Jaka jest ca"lkowita krzywizna okr"egu o promieniu $r$?
Jaka jest ca"lkowita krzywizna krzywej zamkni"etej bez p"etelki,
z jedn"a i z dwoma p"etelkami?
Czy mo"rna uog"olni"c ca"lkowit"a krzywizn"e tak"re na inne krzywe,
ni"r tylko na sk"ladaj"ace si"e z "luk"ow okr"eg"ow?

\subsection{Bardzo trudny problem.}
M"owimy, "re krzywa ma wierzcho"lek w danym punkcie, je"sli krzywizna
ma w tym punkcie lokalne maksimum lub minimum. Jasne jest, "re ka"rda
krzywa zamkni"eta ma co najmniej dwa wierzcho"lki. Poka"r (ale to
naprawd"e trudne!), "re ka"rda zamkni"eta krzywa bez p"etelek
ma co najmniej cztery wierzcho"lki.
%
\section{Geometria zewn"etrzna, geometria~wewn"etrzna.}
Nasza krzywa w poprzednim rozdziale le"ra"la sobie na
p"laszczy"znie. Jej krzywizna zale"ra"la od tego, jak kto"s j"a
tam po"lo"ry"l. Ale czy ta krzywizna mog"laby by"c zmierzona przez
jednowymiarowego mieszka"nca krzywej, kt"ory nic nie wie o niczym, co
wystaje poza t"a krzyw"a? Oczywi"scie nie! Po prostu mo"remy
paluchem poci"agn"a"c za t"a krzyw"a zmieniaj"ac zupe"lnie
jej kszta"lt i w og"ole nie zmieniaj"ac np. d"lugo"sci jednowymiarowych
przedmiot"ow le"r"acych na krzywej. Mo"remy dowolnie zmienia"c krzywizn"e,
a mieszkaniec krzywej nic nie zauwa"ry!

Dlatego nasza krzywizna jest obiektem nale"r"acym do geometrii zew\-n"et\-rznej i
nie nale"r"acym do geometrii wewn"etrznej.
{\bf Geometria zewn"etrzna zajmuje si"e obiektami, kt"ore
zale"r"a od tego, w jaki spos"ob badany "swiat} (w tym wypadku
nasza krzywa) {\bf jest zanurzony w przestrzeni o wi"ekszej liczbie wymiar"ow}
(w naszym przypadku by"la to p"laszczyzna).

Natomiast {\bf geometria wewn"etrzna zajmuje si"e tylko tymi obiektami i w"la\-sno\-"scia\-mi,
 kt"ore
mog"a by"c zbadane i zmierzone przez mieszka"ncow danego "swiata}.
Oczy\-wi"s\-cie takie w"lasno"sci nie zale"r"a od sposobu, w jaki dany
"swiat jest zanurzony w przestrzeni o wi"ekszej liczbie wymiar"ow.

Danym "swiatem mo"re by"c na przyk"lad
krzywa, kt"ora mo"re by"c r"o"rnie po"lo"rona
na p"laszczyznie albo kartka papieru, kt"ora mo"re by"c albo p"lasko
roz"lo"rona albo zwini"eta w kawa"lek bocznej powierzchni walca.
R"o"rne sposoby zanurzenia nie zmieniaj"a oczywi"scie wewn"etrznej
geometrii.

Poniewa"r ka"rda krzywa wygl"ada dla jej mieszka"nca tak samo jak
kawa"lek prostej, wi"ec {\bf krzywa nie ma wewn"etrznej krzywizny}.
B"edziemy wi"ec musieli wybra"c si"e na poszukiwanie krzywizny
wewn"etrznej w wy"rsze wymiary...

\section{Geometria na sferze}
Chcemy uprawia"c geometri"e na sferze. Na p"laszczy"znie mieli"smy
proste i odcinki, czyli najkr"otsze krzywe "l"acz"ace punkty--na
sferze takimi krzywymi s"a okr"egi, ale nie wszystkie. Tylko te,
kt"ore otrzymuje si"e przecinaj"ac sfer"e p"laszczyzn"a
przechodz"ac"a przez "srodek sfery. Takie okr"egi nazywamy
{\bf wielkimi ko"lami}. O jakich s"lynnych wielkich ko"lach uczymy 
si"e na lekcjach geografii? Umawiamy si"e, "re du"rymi literami b"edziemy
oznacza"c punkty, a ma"lymi wielkie ko"la.

Oczywi"scie tak jak na p"laszczy"znie przez ka"rde dwa punkty mo"rna
przeprowadzi"c prost"a, tak {\bf przez ka"rde dwa punkty sfery mo"rna
przeprowadzi"c wielkie ko"lo} (czasami nawet na wiele r"o"rnych
sposob"ow! kiedy?).

O ile jednak na p"laszczy"znie istniej"a proste, kt"ore si"e nie
przecinaj"a (czyli s"a r"ownoleg"le), to {\bf na sferze przecinaj"a
sie ka"rde dwa wielkie ko"la (i to w dw"och punktach!)}.

Na sferze mo"remy mierzy"c odleg"lo"sci r"o"rnych punkt"ow od siebie.
Umawiamy si"e, "re liczymy odleg"lo"s"c chodz"ac po powierzchni
sfery (tak jak ludzie mierz"a odleg"lo"sci na powierzchni Ziemi),
a nie na skr"oty (jak to robi"a krety).

Z mierzeniem k"at"ow miedzy ko"lami wielkimi nie ma chyba problem"ow.

Aha, i jeszcze ma"le za"lo"renie, kt"ore troch"e nam upro"sci rachunki: 
powiedzmy, "re nasza sfera ma promie"n d"lugo"sci $1$.

{\bf Zauwa"r, "re wy"rej rozwa"rane poj"ecia (wielkie ko"lo,
odleg"lo"s"c dw"och punkt"ow, k"at pomi"edzy wielkimi ko"lami)
nale"r"a do geometrii wewn"etrznej sfery.}

\subsection{Co to jest punkt? Jak mierzy"c odleg"lo"s"c punkt"ow?}
\label{punkt}
Co to jest punkt? 
To nie jest trudne pytanie: po wprowadzeniu uk"ladu wsp"o"lrz"ednych
w ten spos"ob, by pocz"atek uk"ladu wsp"o"lrz"ednych by"l
"srodkiem sfery,
mo"remy my"sle"c o punktach sfery jako o tr"ojkach liczb $(x,y,z)$ takich,
"re $x^2+y^2+z^2=1$.

Mamy dane dwa punkty le"r"ace na sferze: $P=(x_1,y_1,z_1)$ i
$Q=(x_2,y_2,z_2)$. Jaka jest ich odleg"lo"s"c?
Zauwa"rmy, "re odleg"lo"s"c dw"och punkt"ow sfery jest
r"ow\-na k"atowi jaki tworz"a te punkty ze "srodkiem sfery.
Z w"lasno"sci iloczynu skalarnego mamy
$$\stackrel{\longrightarrow}{OP}\cdot \stackrel{\longrightarrow}{OQ}=\cos 
\angle POQ \mid OP\mid \ \ \mid OQ \mid$$
Tak wi"ec oznaczaj"ac odleg"lo"s"c naszych punkt"ow przez $l$, mamy
$$\cos l=\stackrel{\longrightarrow}{OP}\cdot 
\stackrel{\longrightarrow}{OQ}=x_1 x_2+y_1 y_2+z_1 z_2$$

\subsection{Co to jest wielkie ko"lo?}
Ka"rde wielkie ko"lo jest wyznaczone przez p"laszczyzn"e, kt"ora przez
nie przechodzi. Og"olne r"ownanie p"laszczyzny le"r"acej w przestrzeni
i przechodz"acej przez "srodek uk"ladu wsp"o"lrz"ednych to
$a x+b y+c z=0$, gdzie $a,b,c$ s"a dowolnymi liczbami
(za wyj"atkiem przypadku, gdy $a,b,c$ s"a wszystkie r"owne $0$).
Oczywi"scie je"sli wszystkie trzy wsp"o"lczynniki przemno"rymy przez
ten sam czynnik, otrzymamy to samo wielkie ko"lo, dla uproszczenia wi"ec
mo"remy przemno"ry"c je tak, "reby $a^2+b^2+c^2=1$.

Tak wi"ec mo"remy my"sle"c o wielkim kole jako o tr"ojce liczb
$(a,b,c)$ takiej, "re $a^2+b^2+c^2=1$.

Przy okazji wida"c, w jaki spos"ob mo"rna "latwo sprawdzi"c, czy prosta
$(a,b,c)$ przechodzi przez punkt $(x,y,z)$: jest tak wtedy, gdy $ax+by+cz=0$.

\subsection{Dualno"s"c.}
Rozpatrzmy tr"ojk"e liczb $(a,b,c)$, $a^2+b^2+c^2=1$. Mo"remy t"e
tr"ojk"e interpretowa"c jako wsp"o"lrz"edne pewnego punktu
sfery $P$ lub jako wsp"o"lczynniki wyznaczaj"ace pewne wielkie ko"lo $p$.
Jaka relacja zachodzi pomi"edzy tak wyznaczonymi punktem i wielkim ko"lem?

To proste: dla dowolnego punktu $Q=(x,y,z)$ le"r"acego na wielkim kole $p$
iloczyn skalarny wektor"ow $\stackrel{\longrightarrow}{OQ}$ i
$\stackrel{\longrightarrow}{OP}$ jest r"owny $ax+by+cz=0$ i te wektory
s"a prostopad"le. Zatem
prosta $OP$ jest prostopad"la do p"laszczyzny wielkiego ko"la $p$.

Wprowadzimy przekszta"lcenie zwane dualno"sci"a: po prostu
punkt o wsp"o"lrz"ednych $(x,y,z)$ zast"apimy wielkim ko"lem
$(x,y,z)$, a wielkie ko"lo $(a,b,c)$ zast"apimy punktem o wsp"o"lrz"ednych
$(a,b,c)$. Punkty przechodz"a  w wielkie ko"la, a wielkie ko"la w punkty.

Udowodnij, "re {\bf je"sli punkt $P$ le"ry na kole $s$ i w wyniku dualno"sci
punkt $P$ przejdzie w ko"lo $p$, a ko"lo $s$ w punkt $S$, to
r"ownie"r punkt $S$ le"ry na kole $p$}.

Udowodnij, "re {\bf je"sli w wyniku dualno"sci punkty $P$, $Q$ przechodz"a
w wielkie ko"la $p$, $q$, to odleg"lo"s"c $PQ$ jest r"owna jednemu z dw"och k"at"ow,
jaki tworz"a przecinaj"ace si"e wielkie ko"la $p$, $q$}.

Dzi"eki dualno"sci mo"remy "latwo obliczy"c k"at $\alpha$,
jaki tworz"a
przecinaj"ace si"e wielkie ko"la $(a_1,b_1,c_1)$ i $(a_2,b_2,c_2)$.
Wiemy, "re ten k"at jest r"owny odleg"lo"sci punkt"ow
$(a_1,b_1,c_1)$ i $(a_2,b_2,c_2)$, a t"e odleg"lo"s"c mierzyli"smy
w rozdziale \ref{punkt} przy pomocy iloczynu skalarnego:
$\cos \alpha=a_1 a_2+b_1 b_2+c_1 c_2$.

Dualno"s"c przeprowadza wierzcho"lki tr"ojk"ata pewne wielkie
ko"la, a boki w pewne punkty. Te wielkie ko"la i punkty tworz"a
tr"ojk"at dualny do pocz"atkowego. {\bf Boki tr"ojk"ata dualnego
s"a r"owne k"atom \underline{przyleg"lym} wyj"sciowego tr"ojk"ata
(k"at przyleg"ly to $\pi$ minus k"at wewn"etrzny) i na 
odwr"ot.}\label{strona}\footnote{{\em Ojej! A dlaczego boki tr"ojk"ata 
dualnego nie s"a po prostu r"owne
k"atom \underline{wewn"etrznym} tr"ojk"ata wyj"sciowego i na odwr"ot?}
Definiuj"ac dualno"s"c prze"slizgn"eli"smy si"e nad do"s"c
wa"rnymi technicznymi szczeg"o"lami. Po prostu jedno ko"lo wielkie
mo"re by"c opisane na {\bf dwa} sposoby jako tr"ojka liczb, np. jako
$(a,b,c)$ lub jako $(-a,-b,-c)$, a zatem przy dualizacji odpowiadaj"a mu
{\bf dwa} punkty. \.Zeby to odpowiednio"s"c by"la wzajemnie
jednoznaczna trzeba na wielkim kole namalowa"c strza"lk"e (nazywamy to
wyborem orientacji). Mo"resz my"sle"c o tym tak, "re na r"owniku Ziemi zaznaczyli"smy jej
kierunek rotacji. Kiedy zastanawiasz si"e, kt"ory z biegun"ow planety
nale"ry wybra"c jako punkt dualny do wielkiego ko"la, powiniene"s zawsze
wybra"c biegun p"o"lnocny (Puchatek go znalaz"l!). Natomiast
gdy rysujesz wielkie ko"lo dualne do jakiego"s punktu, powiniene"s na nim
jeszcze zaznzczy"c strza"lk"e. Zn"ow masz dwie mo"rliwo"sci...
Wybierz wi"ec taki kierunek rotacji, by wyj"sciowy punkt by"l biegunem
p"o"lnocnym!

Czy jeszcze pami"etasz to zadanie:
{\em udowodnij, "re je"sli w wyniku dualno"sci punkty $P$, $Q$ przechodz"a
w wielkie ko"la $p$, $q$, to odleg"lo"s"c $PQ$ jest r"owna {\bf jednemu z dw"och k"at"ow},
jaki tworz"a przecinaj"ace si"e wielkie ko"la $p$, $q$}. Teraz ju"r mo"remy
powiedzie"c, o {\bf kt"ory} k"at chodzi. To tak samo, jak narysowanie strza"lek na
przecinaj"acych si"e prostych czyni z nich p"o"lproste, a k"at mi"edzy
p"o"lprostymi jest jednoznacznie okre"slony.

{\em Dobra, dobra, co to ma wsp"olnego z naszym tr"ojk"atem?}
To proste, boki tr"ojk"ata to kawa"lki wielkich k"o"l. Koniecznie
trzeba namalowa"c na ka"rdym strza"lk"e. Radz"e zrobi"c ci to tak,
by kaza"c jakiej"s mr"oweczce obej"s"c tw"oj tr"ojk"at i kierunki
strza"lek na ka"rdym boku wybra"c zgodnie z kierunkiem, w jakim sz"la
po nim mr"owka. Sprawdz, "re k"aty mi"edzy tak zorientowanymi wielkimi
ko"lami to dok"ladnie k"aty przyleg"le tr"ojk"ata.}

Dualno"s"c okazuje si"e by"c bardzo przydatnym narz"edziem.
Powiedzmy, "re umiemy udowodni"c jakie"s twierdzenie
zachodz"ace dla dowolnego tr"ojk"ata sferycznego.
Je"sli zamienimy
w tym twierdzeniu s"lowa {\em punkt, wierzcho"lek, miara k"ata}
na s"lowa {\em wielkie ko"lo, bok, d"lugo"s"c boku} i
na odwr"ot, a d"lugo"sci bok"ow $a$, $b$, $c$ i miary k"at"ow
$\alpha$, $\beta$, $\gamma$ zast"apimy przez $\pi-\alpha$, $\pi-\beta$, $\pi-\gamma$ i
$\pi-a$, $\pi-b$, $\pi-c$, to dostaniemy inne prawdziwe twierdzenie.
Rozumiesz dlaczego tak jest?
\subsection{Suma k"at"ow w tr"ojk"acie.}
Czy suma k"at"ow w tr"ojk"acie jest r"owna zawsze $\pi$? Od czego
to zale"ry? Do tego problemu jeszcze wr"ocimy w rozdziale
\ref{krzywyswiat}.

\subsection{Twierdzenie Pitagorasa.}
Rysujemy na sferze tr"ojk"at prostok"atny $ABC$. Mo"remy sobie tak
wybra"c uk"lad wsp"o"lrzednych w przestrzeni, "reby $O$--pocz"atek
uk"ladu wsp"o"lrzednych pokrywa"l sie ze "srodkiem sfery. O"s $x$
wybieramy tak, "reby przecina"la sfer"e w punkcie $C$. Osi $y$ i $z$
wybieramy tak, "reby punkt $A$ le"ra"l na p"laszczy"znie wyznaczonej
przez osi $x$ i $y$, a punkt $B$ le"ra"l na p"laszczy"znie wyznaczonej
przez osi $x$ i $z$.

Jakie wsp"o"lrzedne maj"a wierzcho"lki tr"ojk"ata? Oczywi"scie
punkt $C$ ma wsp"o"lrzedne $[1,0,0]$. Odleg"lo"s"c punktu $A$ od $C$,
czyli k"at $\angle AOC$ jest r"owna $b$. Wida"c wi"ec, "re
$A=[\cos b,\sin b,0]$.

Analogicznie dowodzi si"e, "re $B=[\cos a,0,\sin a]$.

Iloczyn skalarny wektor"ow $\stackrel{\longrightarrow}{OA}$ i
$\stackrel{\longrightarrow}{OB}$ jest r"owny kosinusowi d"lugo"sci
boku $c$. Z drugiej jednak strony ten iloczyn jest
r"owny po prostu $\cos a \cos b$. I tak otrzymujemy {\bf twierdznie
Pitagorasa:}
$$ \cos c=\cos a \cos b $$
Czy przypomina ono stare twierdzenie Pitagorasa dla p"laszczyzny?

\subsection{Twierdzenie kosinus"ow (wzory al-Battaniego).}
\label{twkosinusow}
Zapewne znasz twierdzenie kosinus"ow, czyli uog"olnienie twierdzenia Pitagorasa na tr"ojk"aty niekoniecznie prostok"atne.
Spr"obuj samodzielnie odkry"c, jak ono wygl"ada na sferze!

\subsection{Trygonometria.}
\label{trygonometria}
Rozpatrzmy tr"ojk"at prostok"atny $ABC$, $\angle ACB=\frac{\pi}{2}$.
Na p"laszczy"znie mieliby"smy $\sin \alpha=\frac{a}{c}$ oraz $\cos \alpha=\frac{b}{c}$.
A jak jest na sferze?

\subsection{Twierdzenie sinus"ow.}
Udowodnij twierdzenie sinus"ow:{\em w tr"ojk"acie sferycznym o bokach
$a$, $b$, $c$ i k"atach $\alpha$, $\beta$, $\gamma$ zachodzi
$$\frac{\sin \alpha}{\sin a}=\frac{\sin \beta}{\sin b}=\frac{\sin \gamma}{\sin c}$$
}
Bardzo pomys"lowy dow"od znajduje si"e w {\em Podr"eczniku in"ryniera},
ale mo"rna te"r skorzysta"c z trygonometrii.

\subsection{Obw"od i pole ko"la.}
\label{obwod}
Wszyscy znamy wzory na pole ko"la i obw"od okr"egu o danym promieniu.
A jak wygl"adaj"a analogiczne wzory na sferze?\footnote{Mo"re zerkniesz na 
stron"e 90. ksi"a"rki Marka Kordosa, Wyk"lady z historii matematyki,
na kt"orej pokazany jest elementarny spos"ob na liczenie powierzchni sfery.
Mo"resz go potrzebowa"c!}

\subsection{Sfera o dowolnym promieniu.}
Wszystkie wzorki w tym rozdziale dotyczy"ly sfery o promieniu $1$.
Jak nale"ry je zmodyfikowa"c, by dotyczy"ly sfery o dowolnym promieniu
$r$?

\subsection{Ma"ly kawa"lek sfery to prawie p"laszczyzna.}
{\bf Ma"ly kawa"lek sfery jest bardzo podobny do p"laszczyzny, zatem
wszystkie wzory geometrii sferycznej powinny dla ma"lych d"lugo"sci
odcink"ow
sprowadzi"c si"e do zwyk"lych wzor"ow geometrii p"laskiej.}

Sprawd"z to: we"z kilka wzor"ow trygonometrii sferycznej (np. udowodnione
przez nas twierdzenie Pitagorasa dla sfery; inne wzorki mo"resz znale"z"c
w Encyklopedii Szkolnej, w Poradniku in"ryniera i w Poradniku matematycznym).
Co si"e z nimi dzieje, je"sli za"lo"rysz, "re rozmiar
rozwa"ranego tr"ojk"ata jest bardzo ma"ly w por"ownaniu z promieniem
sfery? Pami"etaj, "re dla ma"lych k"at"ow $\alpha$ mamy
$\sin \alpha \approx \alpha-\frac{\alpha^3}{6}$ i $\cos \alpha \approx 
1-\frac{\alpha^2}{2}+\frac{\alpha^4}{24}$.
Ostro"rnie! Ma"ly tr"ojk"at ma ma"le boki, ale ca"lkiem du"re k"aty.

Zamiast zak"lada"c, "re nasz tr"ojk"at jest ma"ly mo"resz za"lo"ry"c, "re 
sfera ma du"ry promie"n.

Patrz"ac cho"cby na twierdzenie Pitagorasa wida"c, "re twierdzenia
geometrii sferycznej wygl"adaj"a do"s"c egzotycznie w por"ownaniu
z ich klasycznymi odpowiednikami. Czasami nawet jedno
twierdzenie o zwyk"lym tr"ojk"acie ma kilka r"o"rnych odpowiednik"ow
dla tr"ojk"ata sferycznego.

\subsection{} Sfera o niesko"nczenie du"rym promieniu to w granicy zwyk"la p"laszczyzna. 
Czy mo"rna jako"s wprowadzi"c dualno"s"c na p"laszczyznie? Nie znam 
odpowiedzi i ch"etnie si"e dowiem.

\subsection{P"laszczyzna rzutowa.}
Definiuj"ac dualno"s"c prze"slizgn"eli"smy si"e nad do"s"c
wa"rnymi technicznymi szczeg"o"lami. Przecie"r istniej"a tak naprawd"e
dwa punkty na sferze, kt"ore s"a dualne do danego wielkiego ko"la.
Te punkty s"a do siebie antypodyczne, czyli le"r"a symetrycznie
wzgl"edem "srodka sfery. Troch"e nie wiadomo, kt"ory wybra"c...

Innym rozwi"azaniem tego problemu ni"r to z przypisu na stronie \pageref{strona}
jest uprawianie geometrii nie na sferze, a na p"laszczyznie rzutowej.
Po prostu umawiamy si"e, "re ka"rde dwa antypodyczne punkty to jeden,
ten sam punkt. Mo"remy o tym my"sle"c nast"epuj"aco:
mamy teraz tylko jedn"a po"lkul"e i je"sli dochodzimy do jej
brzegu i chcemy i"s"c dalej, natychmiast przechodzimy na antypodyczny
punkt.

Sprawd"z, "re p"laszczyzna rzutowa, podobnie jak wst"ega M\"obiusa
jest nieorientowalna.

Mo"rna te"r um"owi"c si"e, "re nie b"ed"a nas interesowa"c
odleg"lo"sci punkt"ow ani miary k"at"ow. {\em No to niby czym
b"edziemy si"e zajmowa"c?} Nie jest a"r tak "zle, przy takich
ograniczeniach mo"remy ci"agle powiedzie"c na przyk"lad, "re
punkt $P$ le"ry na wielkim kole $q$ (zamiast wielkie ko"lo m"owi
si"e te"r prosta), albo "re dane trzy punkty s"a wsp"o"lliniowe,
lub "re trzy proste maj"a punkt wsp"olny. Dzia"l matematyki zajmuj"acy si"e 
takimi problemami nazywa si"e geometri"a rzutow"a. Zawiera ona wiele 
zabawnych twierdze"n,
cho"cby twierdzenie Pappusa albo Pascala. Polecam ksi"a"rk"e Coxetera,
bo z krzywizn"a geometria rzutowa nie ma zbyt du"ro wsp"olnego i nie b"ed"e si"e
ni"a tu zajmowa"c.

\subsection{}
Spr"obuj udowodni"c jakie"s wzory trygonomerii sferycznej, np wz"or na
pole tr"ojk"ata sferycznego prostok"atnego o znanych przyprostok"atnych,
albo pole dowolnego tr"ojk"ata sferycznego przez jego dwa boki i k"at
albo sferyczny odpowiednik wzoru Herona: wyra"z pole tr"ojk"ata sferycznego
przez jego boki. Tyle znasz wzor"ow dla tr"ojk"ata na p"laszczyznie...
Znajd"z ich odpowiedniki! Udowodnij wzory Cagnoliego i Lhuiliera (a co to jest?).

Czy dwusieczne tr"ojk"ata sferycznego si"e przecinaj"a? A symetralne?
"Srodkowe? Wysoko"sci?

Czy potrafisz znale"z"c sferyczne odpowiedniki twierdzenia Menelausa i Cevy?

Tyle twierdze"n, tyle roboty przed tob"a!

%
%
%
%
%
%
\subsection{Mapy.}
Mapa--to odwzorowanie z kawa"lka sfery na p"laszczyzn"e.
To odwzorowanie nie mo"re by"c byle jakie. Ka"rdy, kto widzia"l,
jakie patologiczne funkcje potrafi"a wymy"sla"c matematycy wie, "re
przyda"loby si"e za"lo"ry"c, "re to odwzorowanie jest ci"ag"le.
Ponadto, "reby na mapie nic si"e nie skleja"lo, lepiej za"lo"ry"c,
"re nasze odwzorowanie jest r"o"rnowarto"sciowe.

Mapy sfery mo"rna robi"c na r"o"rne sposoby. Jednym z naj"ladniejszych
jest rzut stereograficzny: sfer"e k"lasziesz na stole na biegunie
po"ludniowym. Ka"rdemu\footnote{Naprawd"e ka"rdemu?} punktowi $P$ sfery
przypisuj"e punkt $P'$, punkt przeci"ecia prostej $NP$ ($N$--to biegun
p"o"lnocny) ze sto"lem. To odwzorowanie przekszta"lca wielkie ko"la
w okr"egi i proste na stole oraz zachowuje k"aty. Jak si"e o tym
przekona"c? Powiem tylko, "re na rzut stereograficzny jest inwersj"a
wzgl"edem pewnej sfery o "srodku $N$. Teraz wystarczy zna"c elementarne
w"lasno"sci inwersji...

Inne ciekawa mapa: sfer"e k"lasziesz na stole na biegunie
po"ludniowym. Prawie ka"rdemu\footnote{Co to znaczy?} punktowi $P$ sfery
przypisuj"e punkt $P'$, punkt przeci"ecia prostej $OP$ ($O$--to "srodek sfery)
ze sto"lem. Niestety, ta mapa nie zachowuje k"at"ow, ale za to
wielkie ko"la przekszta"lca na proste.

Bardzo lubi"e te"r inn"a map"e, kt"ora z kolei zachowuje pole powierzchni 
figur, a otrzymuje si"e j"a rzutuj"ac sfer"e na powierzchni"e walca.

Udowodnij, "re nie istnieje mapa (kawa"lka) sfery, kt"ora zachowywa"laby 
d"lugo"sci krzywych.

Udowodnij, "re mapa dowolnej powierzchni dwuwymiarowej, kt"ora zachowuje 
k"aty i powierzchnie musi te"r zachowywa"c d"lugo"sci krzywych.

\section{P"laszczyzna hiperboliczna.}
Mo"resz sprawdzi"c, "re jaki by"s nie wzi"a"l promie"n sfery,
zawsze g"esto"s"c krzywizny b"edzie dodatnia\footnote{Przepraszam, co to 
jest krzywizna? Ach racja, opowiem o tym w nast"epnych rozdzia"lach.}. A 
szkoda, bo
ch"etnie mia"loby si"e pod r"ek"a jaki"s prosty "swiat o ujemnej 
krzywi"znie. Co to znaczy prosty? To jasne: izotropowy, czyli w ka"rdym 
miejscu wygl"adaj"acy tak samo! Zgodzicie si"e, "re fili"ranka w 
przeciwie"nstwie do sfery ma skomplikowan"a geometri"e w"la"snie przez to, 
"re r"o"rne miejsca wygl"adaj"a r"o"rnie.

Aby skonstruowa"c taki "swiat o ujemnej krzywiznie, przypomnijmy sobie, jak
zrobiona jest sfera.
W przestrzeni tr"ojwymiarowej zadany jest iloczyn skalarny:
$(x_1,y_1,z_1)\cdot(x_2,y_2,z_2)=x_1 x_2+y_1 y_2+z_1 z_2$ dla dowolnych
wektor"ow $(x_1,y_1,z_1)$ i $(x_2,y_2,z_2)$. Je"sli mamy uk"lad 
wsp"o"lrz"ednych,
to mo"remy wed"lug tego samego przepisu oblicza"c iloczyn skalarny
nie tylko wektor"ow, ale i punkt"ow. Sfera jednostkowa to po prostu zbi"or 
takich
punkt"ow $P$, "re $\wektor{OP}\cdot\wektor{OP}=P\cdot P=1$ (przez $O$
oznaczy"lem "srodek uk"ladu wsp"o"lrz"ednych).

Tym razem we"zmy nieco inny iloczyn skalarny: 
$(x_1,y_1,z_1)\cdot(x_2,y_2,z_2)=x_1 x_2+y_1 y_2-z_1 z_2$.
R"o"rnica jest male"nka. Zobacz jednak, jak dziwnie wygl"ada teraz
sfera, zbi"or takich $P$, "re $\wektor{OP}\cdot\wektor{OP}=P\cdot P=-1$. To 
hiperboloida paraboliczna, a my b"edziemy j"a nazywa"c p"laszczyzn"a 
hiperboliczn"a.

Najpierw zobaczmy, jak wygl"adaj"a izometrie, czyli takie przekszta"lcenia
przestrzeni, kt"ore nie zmieniaj"a iloczynu skalarnego.
Kiedy"s (dla przestrzeni ze zwyk"lym iloczynem skalarnym)
by"ly to zwyk"le obroty, np. $(x,y,z) \longmapsto (x \cos \alpha+y \sin 
\alpha,
y \cos \alpha-x \sin \alpha,z)$ i takie przekszta"lcnia s"a r"ownie"r 
izometriami p"laszczyzny hiperbolicznej.

Sprawd"z, "re jedn"a z izometrii dla nowego iloczynu skalarnego jest
$(x,y,z) \longmapsto (x \cosh t+z \sinh t,
y, z \cosh t +x \sinh t)$, gdzie $\cosh t=\frac{e^t+e^{-t}}{2}$.
$\sinh t=\frac{e^t-e^{-t}}{2}$. Jak wygl"ada z"lo"renie takich dw"och
przekszta"lce"n? Czy te"r jest tej postaci?

Zwyk"la sfera wygl"ada w ka"rdym miejscu tak samo, to znaczy sfer"e
mo"rna tak izometrycznie przekszta"lci"c na sfer"e, by dowolny punkt
sfery przekszta"lci"c na dowolny inny. T"e sam"a w"lasno"s"c ma
np. powierzchnia niesko"nczonego walca. Na walcu ka"rdy punkt wygl"ada tak
samo, ale nie ka"rdy kierunek wygl"ada tak samo, bo nie istnieje zbyt wiele izometrii
powierzchni walca, kt"ore zachowywa"ly dany punkt. Powinni"smy wi"ec za"r"ada"c,
by w izotropowym (czyli jednorodnym) "swiecie istnia"la izometria b"ed"aca
obrotem wok"o"l wybranego punktu o dowolny k"at oraz izometria 
przekszta"lcaj"aca dowolny punkt w dowolny inny. Radz"e ci sprawdzi"c, "re 
p"laszczyzna hiperboliczna spe"lnia te wymagania.

Zastan"ow si"e, czym s"a odpowiedniki prostych lub wielkich k"ol
na p"laszczyznie hiperbolicznej.

Zauwa"r, "re nasz nowy iloczyn skalarny nie jest dodatnio okre"slony, to znaczy
czasami $P\cdot P<0$. To smutne, bo to oznacza, "re nie bardzo wiadomo, co
to znaczy $\mid P\mid=\sqrt{P\cdot P}$. Ale je"sli wezmiemy jaki"s wektor
styczny do naszej dziwacznej sfery (powiedzmy, jest to wektor zaczepiony
w punkcie $(0,0,1)$, jest on styczny, wi"ec jest postaci $\vec v=(a,b,0)$),
to "latwo sprawdzi"c, "re wszystko jest w porz"adku, bo $\vec v\cdot
\vec v=a^2+b^2\geq 0$, zatem iloczyn iloczyn skalarny jest dodatni.

Na sferze odleg"lo"s"c $l$ punktu $P$ od $Q$ spe"lnia zale"rno"s"c
$\cos l=\wektor{OP}\cdot\wektor{OQ}$. Bardzo podobnie jest na p"laszczyznie
hiperbolicznej: $\cosh l=\wektor{OP}\cdot\wektor{OQ}$.

Czy pami"etasz, jak robili"smy mapy sfery? Zr"ob dok"ladnie tak samo
mapy p"laszczyzny hiperbolicznej. Mapa analogiczna do rzutu stereograficznego
nazywana jest modelem Poincarego p"laszczyzny hiperbolicznej, a
analogiczna do tej drugiej to model Kleina. Poka"r, "re w modelu Poincarego
proste przechodz"a na proste i okr"egi oraz "re model Poincarego
zachowuje k"aty (zn"ow musisz znale"z"c jak"a"s inwersj"e,
ale uwa"raj! bo nie b"edzie to zwyk"la inwersja), a tak"re "re
w modelu Kleina  proste przechodz"a na proste.
\subsection{Dualno"s"c?}
Hi, hi, hi, a co z dualno"sci"a? Przyznaj"e, "re nie wiem.

\subsection{Magiczna liczba $\sqrt{-1}$.} Popatrz na wzory trygonometryczne 
na sferze i na p"laszczyznie hiperbolicznej. Czy nie s"a do siebie podobne? 
Jak znaj"ac wz"or trygonometrii sferycznej odgadn"a"c odpowiadaj"acy wz"or 
na p"laszczyznie hiperbolicznej?

\section{R"o"rne "swiaty.}
Chcieliby"smy zbada"c geometri"e wewn"etrzn"a powierzchni
dwuwymiarowych (ju"r wiemy, "re geometria wewn"etrzna obiekt"ow jednowymiarowych
jest nieciekawa, a zakrzywione przestrzenie tr"ojwymiarowe mog"a by"c
za trudne jak na pierwsze podej"scie...)

Typowymi przyk"ladami takich powierzchni s"a zwyk"la p"laszczyzna,
bo\-czna powierzchnia walca,
sfera, torus (czyli d"etka), powierzchnia fili"ranki, torus z dwiema dziurkami (ciasteczko z dwoma dziurkami),...

S"a te"r inne, troch"e nieoczekiwane przyk"lady.
Moim ulubionym jest p"la\-szczy\-zna "ruczkowa. Wyobra"zmy sobie
p"laszczyzn"e, na kt"orej mieszka sobie "ruk. Opowiada o tym
Richard Feynmann w doskona"lych {\em Wyk"ladach z fizyki} w tomie II, cz"e"sci 2,
rozdziale 42:
{\it na p"laszczy"znie tej w r"o"rnych miejscach b"edzie
wyst"epowa"la r"o"rna temperatura. Zak"ladamy, "re zar"owno
sam "ruk, jak i pr"ety miernicze, kt"orymi b"edzie si"e on pos"lugiwa"l,
s"a zbudowane z substancji  tym samym wsp"o"lczynniku rozszerzalno"sci
cieplnej. Je"sli nasz "ruk przeniesie pr"et mierniczy w celu wykonania
pomiaru do dowolnego miejsca na p"laszczyznie, pr"et natychmiat, na skutek
rozszerzalno"sci cieplnej, zmieni swoj"a d"lugo"s"c na tak"a, jaka
odpowiada temperaturze w danym miejscu. W podobny spos"ob b"ed"a si"e
zmienia"ly rozmiary ka"rdego obiektu--samego "ruka, pr"eta mierniczego,
tr"ojk"ata czy czegokolwiek innego--gdy"r obiekt taki b"edzie natychmiast
zmienia"l swe rozmiary na skutek rozszerzalno"sci cieplnej. Ka"rdy
przedmiot b"edzie d"lu"rszy w miejscach gor"acych ni"r w miejscach
zimnych, a wszystko b"edzie mia"lo ten sam wsp"o"lczynnik rozszerzalno"sci
cieplnej.}

Mo"remy te"r wyobra"ra"c sobie p"laszcyzn"e "ruczkow"a
troch"e inaczej:
pr"edko"s"c, z jak"a mo"re si"e porusza"c "ruk
zale"ry od miejsca, w kt"orym  si"e znajduje. W zwi"azku z tym znacznie
naturalniejszym dla niego sposobem mierzenia r"o"rnych dr"og na p"laszczyznie
jest nie pomiar ich d"lugo"sci, ale czasu, jaki zabiera "rukowi ich
pokonanie.

Kolejny przyk"lad otrzymujmy tak: bierzemy jak"a"s dwuwymiarow"a 
powierzchni"e, wycinamy jaki"s jej kawa"lek no"ryczkami i w sprytny 
spos"ob sklejamy jej brzeg. W ten spos"ob powstaje znany wielu graczom 
komputerowym torus: dolny brzeg ekranu sklejony jest z g"ornym, a lewy z 
prawym. Podobnym przyk"ladem jest znana nam p"laszczyzna rzutowa.

\section{Geometria wewn"etrzna.}
Zastan"owmy si"e teraz jakie wielko"sci potrafimy mierzy"c
korzystaj"ac wy\-"l"a\-cznie ze "srodk"ow, jakich dostarcza nam geometria
wewn"etrzna powierzchni dwuwymiarowych.

Je"sli interesuje nas geometria wewn"etrzna jakiej"s powierzchni, na
pewno nic nie b"edziemy mogli powiedzie"c o punktach, kt"ore do niej
nie nale"r"a. Podobnie nie jeste"smy w stanie nic powiedzie"c o
wektorach, kt"ore nie s"a styczne do naszej powierzchni.

{\bf Potrafimy zmierzy"c d"lugo"s"c dowolnej krzywej zawartej w
interesuj"acym nas "swiecie.} W szczeg"olno"sci {\bf istnieje najkr"otsza
krzywa "l"acz"aca dwa wybrane punkty--tak"a krzyw"a nazywamy geodezyjn"a}.
I tak, na zwyk"lej p"laszczyznie geodezyjnymi s"a proste, a na sferze
wielkie ko"la.

{\bf Potrafimy te"r zmierzy"c k"at, jaki tworz"a dwie przecinaj"ace si"e
krzywe lub k"at, jaki tworz"a dwa wektory zaczepione w tym samym punkcie.}

A czy umiemy zmierzy"c k"at, jaki tworz"a dwa wektory zaczepione w r"o"r\-nych
punktach? Zobaczmy: je"sli mamy dwa wektory zaczepione w r"o"rnych punktach
na p"laszcyznie, mo"remy zwin"a"c p"laczyzn"e w tr"abk"e tak,
"reby nic nie popsu"c, nic nie zgi"a"c ani nie rozci"agn"a"c.
Oczywi"scie mieszka"ncy p"laszczyzny nic nie zauwa"r"a, a k"at
mi"edzy wektorami si"e zmieni. Tak wi"ec {\bf k"at pomi"edzy
dwoma wektorami zaczepionymi w r"o"rnych punktach nie nale"ry do geometrii
wewn"etrznej.}

Tak samo "latwo jest si"e przekona"c, "re {\bf mo"remy doda"c
dwa wektory, tylko wtedy, gdy s"a zaczepione w tym samym punkcie.}
%

\subsection{Ekwidystanty}
Je"sli mamy jak"a"s krzyw"a, mo"remy znale"z"c zbi"or punkt"ow, kt"ore 
znajduj"a si"e w ustalonej odlg"lo"sci od tej krzywej. Zazwyczaj ten zbi"or 
jest te"r jak"a"s krzyw"a i nazywamy j"a ekwidystant"a (czyli po polsku 
r"ownoodlg"l"a).

Na p"laszczyznie ekwidystant"ami prostej s"a proste do niej r"ownoleg"le, 
a na sferze ekwidystant"ami wielkiego ko"la---r"ownika' s"a r"ownole"rniki.

Co jest ekwidystant"a ekwidystanty?

Co jest ekwidystant"a punktu?

Udowodnij, "re zawsze ekwidystanty danej krzywej tworz"a k"at prosty z 
geodezyjnymi prostopad"lymi do tej krzywej.

\subsection{} Udowodnij, "re na torusie (nawet troch"e pognicionym) istnieje 
co najmniej 5 geodezyjnych, kt"ore s"a zamkni"etymi krzywymi.

\section{Krzywy "swiat.}
\label{krzywyswiat}
\subsection{K"at obrotu.}
Spacerujesz na p"laszczyznie po brzegu wielok"ata. Sw"oj nos
zawsze trzymasz w kierunku marszu, wi"ec w ka"rdym wierzcho"lku
musisz si"e nieco obr"oci"c.
Po powrocie do
punktu startu tw"oj ca"lkowity obr"ot, czyli suma k"at"ow,
o jakie si"e obr"oci"le"s wynosi $2 \pi$.
Na zakrzywionych powierzchniach zazwyczaj ju"r tak nie jest.

\subsection{Suma k"at"ow tr"ojk"ata. Defekt.}
Wszyscy wiemy, "re na p"laszczy"znie suma k"at"ow tr"ojk"ata
jest r"owna $\pi$. Na zakrzywionych powierzchniach zazwyczaj ju"r tak nie
jest
, wi"ec {\bf wygodnie jest
zdefiniowa"c sobie defekt tr"ojk"ata jako sum"e jego k"at"ow
odj"a"c $\pi$}.

Okazuje si"e, "re {\bf defekt tr"ojk"ata na sferze
jest zawsze dodatni i r"owny polu tego tr"ojk"ata} (przypominam, "re nasza sfera ma promie"n $1$).
Bardzo "ladnie jest to pokazane w artykule z {\it Delty}.

Wida"c wi"ec, "re {\bf na sferze suma k"at"ow
tr"ojk"ata jest zawsze wi"eksza od $\pi$}.

\subsection{Defekt wielok"at"ow.}
Zastan"ow si"e, jak zdefiniowa"c defekt dla innych wielok"at"ow.
Udowodnij, "re {\bf defekt dowolnego wielok"ata na sferze o promieniu $1$
jest zawsze dodatni i r"owny polu tego wielok"ata}.

\subsection{Przeniesienie r"ownoleg"le.}
Na p"laszczyznie mamy dany pewien wektor zaczepiony w jakim"s punkcie i
naszym celem jest przeniesienie go do jakiego"s innego punktu. Nie jest
to trudne: rysujemy odcinek pomi"edzy tymi dwoma punktami i
w punkcie docelowym rysujemy wektor o tej samej d"lugo"sci i tworz"acy
ten sam k"at z odcinkiem, co nasz wyj"sciowy wektor.

Mo"remy zrobi"c to samo w troch"e bardziej skomplikowany spos"ob:
rysujemy "laman"a "l"acz"ac"a punkt pocz"atkowy i ko"ncowy
i nasz wektor przesuwamy po kolei przez wszystkie odcinki "lamanej.
T"e operacj"e nazywamy {\bf przesuni"eciem r"ownoleg"lym,
przesun"eli"smy r"ownolegle wektor wzd"lu"r "lamanej}.

Na zakrzywionej powierzchni mo"remy robi"c to samo, rol"e odcink"ow
przejm"a wtedy geodezyjne. Sprawd"z, "re o ile na p"laszczyznie
wynik nie zale"ra"l od tego, po jakiej "lamanej przesuwali"smy wektor,
to na sferze ju"r tak nie jest.

T"e operacj"e nazywamy przeniesieniem r"ownoleg"lym, chocia"r
cho"cby na sferze je"sli przeniesiemy jaki"s wektor gdzie indziej, zazwyczaj
pocz"atkowy i ko"ncowy wektor nie s"a r"ownoleg"le!


\subsection{Defekt, k"at obrotu i przeniesienie r"ownoleg"le.}
Okazuje si"e "re defekt, k"at obrotu i przeniesienie r"ownoleg"le
s"a bardzo ze sob"a silnie zwi"azane.
Rozwa"rmy mianowicie pewien wielok"at na zakrzywionej powierzchni.
Obliczmy jego defekt. Przespacerujmy si"e wzd"lu"r jego ca"lego
obwodu i obliczmy ca"lkowity k"at obrotu odj"a"c $2\pi$.
Przesu"nmy jaki"s wektor wzd"lu"r jego obwodu i zmierzmy k"at mi"edzy
wyj"sciowym i tak przesuni"etym wektorem. Udowodnij, "re te
trzy liczby s"a r"owne.
Je"sli wi"ec b"edziesz  co"s chcia"l udowodni"c o defekcie,
by"c mo"re tw"oj dow"od b"edzie "latwiejszy, je"sli r"ownowa"rnie
udowodnisz podobn"a w"lasno"s"c przesuni"ecia r"ownoleg"lego lub
ca"lkowitego k"ata obrotu.

Zauwa"r, "re przesuni"ecie r"ownoleg"le, defekt i k"at obrotu
nale"r"a do geometrii wewn"etrznej.

\subsection{Defekt jest addytywny.}
Wielko"s"c addytywna to taka wielko"s"c, kt"ora obliczona dla 'sumy'
dw"och obiekt"ow jest r"owna sumie wielko"sci obliczonych dla sk"ladnik"ow,
o ile tylko te sk"ladniki s"a 'roz"l"aczne'. Pole powierzchni,
obj"eto"s"c, masa, "ladunek elektryczny s"a przyk"ladami takich
wielko"sci.

Udowodnij, "re {\bf defekt jest addytywny}, czyli suma defekt"ow dw"och
wielok"at"ow stykaj"acych si"e tylko wzd"lu"r bok"ow jest
r"owny defektowi wielok"ata b"ed"acego sum"a tamtych.

\subsection{Krzywizna wewn"etrzna.}
Widzimy, "re defekt ma wprost cudowne w"lasno"sci: na p"laszczyznie jest stale
r"owny zero i ma niezerowe warto"sci na sferze, zatem jest jako"s
zwi"azany
z tym, co intuicyjnie "l"aczymy ze s"lowem 'krzywizna'.
Ponadto defekt wielok"ata mo"re by"c zmierzony przez mieszka"nc"ow
danej powierzchni, zatem nale"ry do geometrii wewn"etrznej.
Co wi"ecej, defekt jest addytywny, zatem du"re obszary
maj"a zazwyczaj wi"ekszy defekt, ni"r ma"le, co zgadza si"e
z nasz"a intuicj"a, "re efekty zwi"azane z krzywizn"a
staj"a si"e "latwiej obserwowalne na du"rych obszarach.
Mo"remy wi"ec pomy"sle"c, "re defekt mierzy nam ile 'krzywizny'
znajduje si"e wewn"atrz badanego wielok"ata, tak jakby
krzywizna by"la r"ownie namacaln"a wielko"sci"a jak masa czy "ladunek
elektryczny. Dlatego te"r zdefiniujmy {\bf ca"lkowit"a krzywizn"e
danego obszaru powierzchni jako jego defekt}.

Tak jak stosunek masy do obj"eto"sci jest nazywany g"esto"sci"a,
a stosunek "ladunku elektrycznego do obj"eto"sci g"esto"sci"a
"ladunku, tak te"r {\bf stosunek defektu do pola powierzchni nazwijmy
g"esto"sci"a krzywizny lub po prostu krzywizn"a} (uwa"raj, "reby
nie pomyli"c poj"ecia ca"lkowitej krzywizny jakiego"s obszaru i
krzywizny tzn. g"esto"sci krzywizny!).
Aby zmierzy"c g"esto"s"c krzywizny powierzchni dwuwymiarowej w jakim"s
punkcie,
narysujemy ma"ly wielok"at zawieraj"acy nasz punkt. Obliczamy defekt
wielok"ata i dzielimy przez jego pole.

Wielk"a zalet"a powy"rszej definicji jest to, "re tak zdefiniowana
krzywizna jest oczywi"scie w"lasno"sci"a wewn"etrzn"a powerzchni.
Problem polega na tym, "re korzystaj"ac z niej trudno jest zazwyczaj
porachowa"c krzywizn"e jakich"s konkretnych powierzchni (sfera i
p"laszczyzna s"a wyj"atkami--dla nich obliczenie krzywizny nie
powinno by"c "radnym problemem). To dlatego, "re geodezyjne s"a
opisywane przez do"s"c paskudne r"ownania.

\subsection{Krzywizna konkretnych powierzchni.}
Oblicz g"esto"s"c krzywizny na p"laszczyznie i na walcu.
Oblicz g"esto"s"c krzywizny sfery o promieniu $r$ i jej ca"lkowit"a
krzywizn"e.
Oblicz g"esto"s"c krzywizny bocznej powierzchni sto"rka, 
'torusa graczy komputerowych', p"laszczyzny rzutowej, p"laszczyzny 
hiperbolicznej.

\subsection{}
Udowodnij, "re je"sli
krzywizna jakiego"s "swiata wsz"edzie znika, to prostok"at (czyli figura 
z"lo"rona z czterech geodezyjnych, ma wszystkie k"aty proste) ma r"owne 
przeciwleg"le boki. 

Udowodnij, "re taki "swiat jest lokalnie izometryczny z p"laszczyzn"a.

\subsection{A jednak si"e krzywi!} Firma geodezyjna P.Swarthy and 
S.R.French wynajmuje do prac nad swoimi mapami stypendyst"ow Funduszu, 
kt"orzy nigdy nie s"lyszeli o tym, "re Ziemia nie jest p"laska.

Oszacuj, jakiej wielko"sci deformacje b"ed"a zawiera"c mapy obszar"ow o 
wielko"sci: $10$ km, $100 km$, $1000 km$, $10000$ km.

Innymi s"lowy spr"obuj narysowa"c na sferze co"s, co przypomina prostok"at 
i oszacuj na ile jego k"aty nie s"a proste i na ile przeciwlg"le boki nie 
s"a r"owne. 

\subsection{Spacery.}
\label{spacery}
Przespaceruj si"e jako"s po p"laszczy"znie tak, "reby dok"ladnie
wr"oci"c do punku startu i jeszcze "reby by"c mie"c nos skierowany
tak, jak na pocz"atku.

A teraz przenie"s si"e na sfer"e i wykonuj tyle samo krok"ow i te same obroty.
Tym razem jednak prawie na pewno nie wr"oci"le"s do
punktu wyj"scia. I nos jest troch"e obr"ocony...
Mo"remy oczekiwa"c, "re im kr"otsze kroki b"edziesz robi"c, tym te
efekty b"ed"a mniejsze.

Spr"obuj zbada"c to troch"e dok"ladniej: jak (w przybli"reniu) odleg"lo"s"c
po\-cz"a\-tko\-we\-go i ko"ncowego po"lo"renia zale"r"a od d"lugo"sci kroku?
Jak od d"lugo"sci kroku zale"ry k"at, o jaki obraca si"e nasz nos?

{\em Wskaz"owka: wybierz sobie jaki"s konkretny, "latwy do oblicze"n
kszta"lt spaceru. Naj"latwiejsze rachunki otrzymuje si"e, je"sli
spacer ma kszta"lt okr"egu lub prostok"ata.}

Czy umiesz przeliczy"c to samo (oczywi"scie w przybli"reniu) na jakich"s
innych zakrzywionych "swiatach?

\subsection{}
Udowodnij, "re je"sli defekt jest proporcjonalny do pola, to nasz "swiat
jest lokalnie izometryczny z p"laszczyzn"a, sfer"a lub p"laszczyzn"a
hiperboliczn"a.

\section{Promie"n krzywizny krzywej le"r"acej na powierzchni.}
\label{promien}
Rozwa"rmy jak"a"s krzyw"a le"r"ac"a na naszej powierzchni.
Jak zmierzy"c promie"n krzywizny w jakim"s punkcie?

Popatrzmy na kilka przyk"lad"ow. Na bocznej powierzchni walca narysujmy
okr"ag przecinaj"ac walec p"laszczyzn"a prostopad"l"a do
jego osi. Promie"n tego okr"egu jest r"owny promieniowi walca. Ale to
wcale nie jest promie"n krzywizny tego okr"egu zmierzony przez mieszka"nca
powierzchni walca! Przecie"r dla niego ten okr"ag jest prost"a, zatem
wyznaczy on krzywizn"e tego okr"egu jako r"own"a zeru.
Wida"c wi"ec, "re chc"ac zdefiniowa"c promie"n krzywizny krzywej
le"r"acej na pewnej powierzchni tak, by ten promie"n m"og"l by"c
wyznaczony przez geometri"e wewn"etrzn"a, musimy by"c ostro"rni.

Ja proponuj"e nast"epuj"ac"a metod"e mierzenia promienia krzywizny:
je"sli mamy jak"a"s krzyw"a i punkt, w  kt"orym interesuje
nas promie"n krzywizny, przespacerujmy si"e wzd"lu"r kr"otkiego
kawa"lka krzywej otaczaj"acego wybrany przez nas punkt. Po spacerze
mierzymy d"lugo"s"c drogi, po jakiej porusza"la si"e nasza lewa
(oznaczmy t"a drog"e przez $l$)
i prawa stopa (oznaczmy t"a drog"e przez $p$). Je"sli rozstaw st"op jest
r"owny $d$, to promie"n krzywizny jest r"owny
$R=\frac{l d}{l-p}$.

Proponuj"e, "reby"s sprawdzi"l, czy ta definicja daje na p"laszczyznie
dobre wyniki.

Innym dobrym "cwiczeniem b"edzie sprawdzenie, jak dzia"la nasza
definicja na sferze. Rozwa"r okr"ag o promieniu $a$ (promie"n zosta"l
zmierzony przez mieszka"nc"ow sfery). Oblicz promie"n krzywizny tego okr"egu.

Udowodnij wa"rn"a w"lasno"s"c geodezyjnej, "re mianowicie jej
krzywizna
jest
r"owna $0$.

\section{P"laszczyzna "ruczkowa}
\subsection{Dlaczego id"ac prosto idziemy tak krzywo?}
Rozka"r "rukowi i"s"c prosto! Dra"n nie chce. To proste: "ruk zawsze
skr"eca, gdy jego lewe nogi s"a na obszarze, w kt"orym mog"a rozwin"a"c
inn"a pr"edko"s"c, ni"r prawe nogi.

Napisz r"ownanie opisuj"ace zmian"e wektora pr"edko"sci "ruka w 
czasie.\footnote{Zauwa"r, "re s"a dwa powody, dla kt"orych wektor 
pr"edko"sci
mia"lby si"e zmieni"c: "ruk w"lazi na obszar, gdzie mo"re szybciej
lub wolniej biec wi"ec wektor pr"edko"sci si"e wyd"lu"ra lub
skraca oraz "ruk, a z nim i wektor pr"edko"sci, skr"eca z powodu
innej pr"edko"sci lewych i prawych n"og.} Wynik zapisz w postaci:
$$\frac{d}{dt}v^i(t)+\sum_i\sum_j\Gamma^i_{j k} v^j v^k=0,$$
gdzie $v^i$ dla $i\in\{ 1,2\}$ to sk"ladowe wektora $\vec{v}$, w sumach wska"zniki
$j,k$ przebiegaj"a zbi"or indeksuj"acy wsp"o"lrz"edne $\{ 1,2\}$.
$\Gamma$ powinny by"c tylko funkcjami zale"rnymi od funkcji opisuj"acej
szybko"s"c "r"o"lwia i jej pochodnych. Fachowcy nazywaj"a $\Gamma$
symbolami Christofela.

Rozwa"r p"laszczyzn"e "ruczkow"a, na kt"orej pr"edko"s"c
"ruczka zale"ry tylko od jednej wsp"o"lrz"ednej. Napisz r"ownanie,
jakie spe"lnia geodezyjna. Je"sli zamiast robaczka wezmiemy "swiat"lo,
powiniene"s uzyska"c prawo Snella znane ci z fizyki, a m"owi"ace
o stosunku sinus"ow k"ata padania i za"lamania "swiat"la na
granicy o"srodk"ow.

\subsection{Jak mi"lo jest by"c "rukiem.}
Istnieje pi"ekne twierdzenie m"owi"ace o tym, "re ka"rdy "swiat
(dok"ladniej: ka"rda rozmaito"s"c dwuwymiarowa) jest lokalnie
izometryczny  z pewn"a p"laszczyzn"a "ruczkow"a (czyli wok"o"l ka"rdego punktu mo"rna wyci"a"c
kawa"lek rozmaito"sci tak, by by"l konforemnie izometryczny z p"laszczyzn"a).
Nic dziwnego, "re "ruki (dok"ladniej: owady) osi"agn"e"ly taki sukces ewolucyjny.

\section{Charakterystyka Eulera--Poincarego. Wz"or Gaussa--Boneta.}
Rozpatrzmy dowoln"a dwuwymiarow"a powierzchni"e bez brzegu,
na przyk"lad sfer"e, torus, p"laszyzn"e rzutow"a, torus
z kilkoma uszami...
Dzielimy nasz"a powierzchni"e
w dowolny spos"ob na $S$ wielok"at"ow, kt"orych bokami s"a geodezyjne.
Przez $K$ oznaczmy ilo"s"c wszystkich kraw"edzi, a przez $W$ ilo"s"c
wszystkich wierzcho"lk"ow. Krzywizna ca"lkowita naszej powierzchni to
po prostu suma defekt"ow wszystkich wielok"at"ow.
Defekt $n-$k"ata jest r"owny z definicji sumie jego k"at"ow wewn"etrznych
odj"a"c $(n-2)\pi=n \pi-2 \pi$. Suma wszystkich k"at"ow
wychodz"acych z danego wierzcho"lka jest r"owna $2 \pi$, zatem suma
k"at"ow wewn"etrznych wszystkich wielok"at"ow jest r"owna $2 \pi K$.
Suma liczby bok"ow wszystkich wielok"at"ow jest r"owna $2 K$, gdy"r
ka"rda kraw"ed"z nale"ry do dw"och wielok"at"ow. Wida"c wi"ec,
"re suma defekt"ow wszystkich wielok"at"ow jest r"owna
$2 \pi (W-K+S)$.

T"e wielko"s"c, $W-K+S$ oznaczamy tradycyjnie liter"a greck"a $\chi$
(chi) i nazywamy {\bf charakterystyk"a Eulera--Poincarego}. Tak wi"ec nasz rezultat
mo"remy zapisa"c w postaci
$$ \mbox{krzywizna ca"lkowita}=2 \pi \chi.$$
Jest to s"lynny wz"or Gaussa--Boneta.

Udowodnij, "re charakterystyka Eulera nie zale"ry od sposobu, w jaki podzielili"smy
na pocz"atku nasz"a powierzchni"e na wielok"aty.

Oblicz charakterystyk"e Eulera sfery, torusa, torusa z dwoma uchami,
p"laszczyzny rzutowej.

Udowodnij, "re je"sli lekko zmodyfikujemy jak"a"s powierzchni"e
(wgniatanie jakiej"s cz"e"sci paluchem itp.) to nie zmieni to jej
ca"lkowitej krzywizny.

\section{Co by by"lo...}
Co by by"lo, gdyby Ziemia nie spoczywa"la na grzbietach czterech "r"o"lwi,
ale by"la niesko"nczon"a p"laszczyzn"a, p"laszczyzn"a hiperboliczn"a
lub (o zgrozo!) sfer"a? Rozwa"r ekonomiczne, polityczne, historyczne nast"epstwa.

Co by by"lo, gdyby Wszech"swiat mia"l kszta"lt niesko"nczonej p"laskiej przestrzeni,
przestrzeni hiperbolicznej lub tr"ojwymiarowej sfery? Jak wygl"ada"lyby
prawa fiyzki: prawo Coulomba, prawo ci"a"renia...

Zerknij do ksi"a"reczki Gamowa. On te"r si"e nad tym zastanawia.
Czy zgadzasz si"e z wszystkimi jego wnioskami?

\section{Powierzchnie zawarte w przestrzeni tr"ojwymiarowej}
\subsection{Geodezyjna i przeniesienie na powierzchniach zanurzonych w 
przestrzeni tr"ojwymiarowej.}
\label{trojwymiar}
Rozwa"rmy powierzchni"e zanurzon"a w przestrzeni tr"ojwmiarowej.
Narysujmy na niej jak"a"s geodezyjn"a i pozaznaczajmy w ka"rdym
jej punkcie wektor styczny do geodezyjnej. Oczywi"scie te wektory nie musz"a by"c
r"ownoleg"le do siebie; przebiegaj"ac geodezyjn"a wektor styczny
b"edzie si"e jako"s obraca"l. Okazuje si"e, "re {\bf pochodna
z wektora stycznego (czyli po prostu wektor, w kierunku kt"orego obraca si"e
wektor styczny) jest prostopad"ly do naszej powierzchni}. Sprawd"z, "re istotnie jest tak na sferze i na bocznej powierzchni walca.
Teraz udowodnimy to w og"olnym przypadku.

Przypu"s"cmy bowiem, "re tak nie jest, to znaczy pochodna z wektora stycznego
ma sk"ladow"a r"ownoleg"l"a do powierzchni. W takim razie poci"agnijmy
troszeczk"e "srodkowy kawa"lek naszej geodezyjnej w kierunku wyznaczon"a przez t"a
sk"ladow"a. Otrzymamy w ten spos"ob krzyw"a nadal le"r"ac"a na naszej
powierzchni, "l"acz"ac"a te
same punkty ko"ncowe, ale kr"otsz"a od wyj"sciowej geodezyjnej (dlaczego?
popatrz na krzyw"a le"r"ac"a na p"laszczyznie i wykonaj dla niej
t"e sam"a operacj"e). Otrzymali"smy sprzeczno"s"c.

Wiedz"ac to, udowodnij, "re {\bf je"sli przenosimy r"ownolegle jaki"s wektor wzd"lu"r
geodezyjnej, zaznaczymy w ka"rdym punkcie geodezyjnej tak otrzymany wektor
i b"edziemy przebiega"c ca"l"a geodezyjn"a i patrze"c, jak
obraca si"e nasz wektor, to pochodna z tego wektora jest prostopad"la
do naszej powierzchni}. 

\subsection{Krzywizna zewn"etrzna powierzchni.}
\label{zewnetrzna}
Rozwa"rmy jak"a"s dwuwymiarow"a powierzchni"e zawart"a w
zwyk"lej przestrzeni tr"ojwymiarowej. Oczywi"scie krzywizna wewn"etrzna
nie opisuje w pe"lni sposobu, w jaki powierzchnia jest po"lo"rona w
przestrzeni, podobnie jak geometria wewn"etrzna krzywej le"r"acej
na powierzchni nie opisuje jej promienia krzywizny. Chcemy wi"ec tak
zdefiniowa"c krzywizn"e zewn"etrzn"a powierzchni, by zawrze"c
w niej ca"l"a istotn"a informacj"e. Ale jak to zrobi"c?

Na szcz"e"scie jeste"smy prawdziwymi ekspertami od krzywizny krzywej
le"r"acej na p"laszczyznie. Gdyby"smy tak problem
wyznaczenia krzywizny powierzchni potrafili jako"s sprowadzi"c do
mierzenia krzywizny krzywych...

B"edziemy przecina"c nasz"a powierzchni"e p"laszczyzn"a
przechodz"ac"a przez wybrany punkt i prostopad"l"a do powierzchni.
Patrz"ac na p"laszczyzn"e widzimy kawa"lek badanej
powierzchni--krzyw"a b"ed"ac"a cz"e"sci"a p"laszczyzny
i powierzchni. Krzywizn"e tej krzywej oczywi"scie mo"remy policzy"c.
W trakcie przekr"ecania p"laszczyzny ta krzywizna b"edzie si"e
zmienia"c...

Wygodnie
b"edzie wybra"c taki uk"lad wsp"o"lrz"ednych, "reby nasz
punkt znalaz"l si"e w jego "srodku $O$ oraz "reby nasza powierzchnia
by"la styczna do p"laszczyzny $OXY$. Powiedzmy, "re w tym uk"ladzie
wsp"o"lrz"ednych nasza powierzchnia jest opisana r"ownaniem $z=f(x,y)$.
\.Zeby by"lo jeszcze "latwiej, niech $f(x,y)=a x^2+2b xy+c y^2$
\footnote{Je"sli jeste"s pedantem, przelicz to w og"olnym przypadku zak"ladaj"ac tylko,
"re $f$ jest dwukrotnie r"o"rniczkowalna i "re jej pierwsze pochodne
znikaj"a w zerze. Dostaniesz to samo.}.

Teraz b"edziemy przecina"c nasz"a powierzchni"e p"laszczyznami
prze\-cho\-dz"a\-cy\-mi przez punkt $O$ i prostopad"lymi do p"laszczyzny $OXY$.
Na takiej p"la\-szczy\-znie wygodnie jest mie"c uk"lad wsp"o"lrz"ednych.
Jedn"a o"s ju"r mamy: jest ni"a o"s $OZ$. Jako drug"a o"s $OT$
bierzemy prost"a b"ed"ac"a przeci"eciem naszej p"laszczyzny
i p"laszczyzny $OXY$. Pomi"edzy wsp"o"lrz"ednymi $(x,y,z)$ i
$(t,z)$ punktu le"r"acego na naszej p"laszcyznie jest prosty zwi"azek:
$x=\cos \alpha t$, $y=\sin \alpha t$, $z=z$ (przez $\alpha$ oznaczyli"smy
k"at $\angle XOT$).

Patrz"ac na nasz"a p"laszczyzn"e widzimy kawa"lek badanej
powierzchni--krzy\-w"a b"ed"ac"a cz"e"sci"a p"laszczyzny
i powierzchni. R"ownanie tej krzywej to $z(t)=f(\cos \alpha t,\sin \alpha t)=
(a \cos^2\!\alpha+b \sin \alpha \cos \alpha+c \sin^2\!\alpha) t^2$. Ci, 
kt"orzy policzyli
krzywizn"e paraboli wiedz"a, "re krzywizna naszej krzywej jest r"owna
$2 (a \cos^2\!\alpha+b \sin \alpha \cos \alpha+c \sin^2\!\alpha)$.
%
Krzywizna mierzona w pewnych kierunkach jest wi"eksza, w innych mniejsza.
Udowodnij, "re
kierunek, w kt"orym krzywizna jest najwi"eksza, jest prostopad"ly
do kierunku, w kt"orym krzywizna jest najmniejsza. Te kierunki nazywamy
{\bf kierunkami g"l"ownymi}. Krzywizny w kierunkach g"l"ownych nazywamy
{\bf krzywiznami g"l"ownymi}.

Oblicz g"l"owne krzywizny dla p"laszczyzny, sfery, powierzchni bocznej
walca.

Oblicz g"l"owne krzywizny dla rozwa"ranej powierzchni $f(x,y)=a x^2+b xy+c 
y^2$
w punkcie $(0,0,0)$.

Udowodnij, "re znaj"ac kierunki g\l"owne i krzywizny g"l"owne
jeste"smy w stanie obliczy"c krzywizn"e w dowolnym kierunku.

\subsection{A gdzie ta prawdziwa krzywizna?}
A gdzie ta prawdziwa krzywizna? To zale"ry, kto pyta.

Badacza baniek mydlanych interesuje tzw. krzywizna "srednia, czyli
$k_{min}+k_{max}$ (dlaczego?)

A krzywizna wewn"etrzna? Jak j"a obliczy"c z $k_{min}$ i $k_{max}$?
To tajemnica, obliczymy to troch"e p"o"zniej.

\subsection{Odwzorowanie Gaussa.}
%
Dwuwymiarowa powierzchnia zanurzona w przestrzeni tr"ojwymiaroej ma
dwie strony\footnote{A co si"e dzieje, gdy nasza powierzchnia ma tylko 
jedn"a
stron"e, podobnie jak wst"ega M\"obiusa? Nie denerwujemy si"e,
poniewa"r interesuj"a nas tylko lokalne w"lasno"sci powierzchni,
mo"remy no"ryczkami wyci"a"c ma"ly kawa"lek, kt"ory b"edzie
mia"l dwie strony.}. Je"sli rysujemy jednostkowy wektor normalny do naszej
 powierzchni,
mo"remy to zrobi"c na dwa sposoby, wybieramy wi"ec jedn"a ze stron
i w t"e stron"e b"edziemy wybiera"c zwrot wektor"ow normalnych.

Odwzorowanie Gaussa przyporz"adkowuje punktom powierzchni punkty sfery
w nast"epuj"acy spos"ob: je"sli mamy jaki"s punkt powierzchni,
wyznaczamy wektor, kt"ory jest normalny do powierzchni w tym w"la"snie
punkcie. Nast"epnie ten wektor przesuwamy tak, by jego punkt zaczepienia
znalaz"l si"e w "srodku naszej sfery. Teraz koniec wektora wskazuje
ten punkt sfery, o kt"ory nam chodzi"lo.

Analogicznie mo"remy odwzorowanie Gaussa te"r zdefiniowa"c dla wektor"ow
stycznych do naszej powierzchni tak, by w analogiczny spos"ob przypisywa"c
im wektory styczne do sfery. Po prostu patrzymy na punkt zaczepienia wektora
stycznego do naszej powierzchni, przekszta"lcamy ten punkt zwyk"lym przekszta"lceniem
Gaussa i tak dostajemy le"r"acy na sferze punkt zaczepienia nowego wektora.
W nowym punkcie zaczepienia rysujemy wektor r"ownoleg"ly do wyj"sciowego--i ju"r!

Rozwa"rmy dowolny obszar na naszej powierzchni, kt"orego brzeg jest jak"a"s
krzyw"a. We"zmy jaki"s wektor le"r"acy na
tej krzywej.
Wiemy, "re je"sli przesuniemy go wzd"lu"r tej krzywej, to k"at, o jaki
ten wektor si"e obr"oci jest r"owny defektowi naszego obszaru.

Odwzorowanie Gaussa przeprowadza nasz obszar w jaki"s kawa"lek sfery,
krzywa ograniczaj"aca nasz obszar zastaje przeprowadzona w krzyw"a
ograniczaj"ac"a ten kawa"lek sfery. R"ownie"r wektory le"r"ace
wzd"lu"r krzywej zostaj"a przekszta"lcone na wektory zaczepione na sferze.

Niezwyk"le jest to, "re zgodnie z tym, co pokazali"smy w rozdziale \ref{trojwymiar},
te wektory na sferze s"a r"ownie"r przenoszone r"ownolegle, zatem
defekty obu obszar"ow s"a takie same. Innymi s"lowy: {\bf odwzorowanie
Gaussa zachowuje krzywizn"e wewn"etrzn"a}. Poniewa"r jednak na
sferze defekt dowolnego obszaru jest r"owny jego powierzchni, zatem
{\bf odwzorowanie Gaussa przeprowadza dowolny obszar powierzchni w pewien
obszar le"r"acy na sferze, kt"orego pole jest r"owne defektowi
wyj"sciowego obszaru}. Zatem {\bf g"esto"s"c krzywizny danej
powierzchni jest r"owna czynnikowi, o jaki powi"eksza si"e pole powierzchni
przy odwzorowaniu Gaussa}.

\subsection{Odwzorowanie Gaussa w dzia"laniu.}
Rozpatrzmy ponownie powierzchni"e $z=a x^2+2bxy+c y^2$.
\L{}atwo sprawdzi"c, "re wektorem normalnym do naszej powierzchni jest
wektor $(-2ax-2by,-2bx-2cy,1)$. Co prawda nie ma on zazwyczaj d"lugo"sci
$1$, ale nie b"edziemy si"e tym przejmowa"c, bo dla otoczenia
punktu $(0,0,0)$ ta r"o"rnica nie
jest istotna\footnote{Oczywi"scie mo"rna dalsze rachunki przeprowadzi"c dla
tego wektora ju"r po przeskalowaniu go tak, by mia"l d"lugo"s"c $1$,
ale rachunki b"ed"a troch"e gorsze... Mo"resz to sprawdzi"c.}.
Ten wektor wyznacza nam odwzorowanie Gaussa. Mo"remy o nim my"sle"c jak
o odwzorowaniu kawa"lka p"laszczyzny (bo nasza powierzchnia $z=f(x,y)$
lokalnie bardzo przypomina p"laszczyzn"e $z=0$) w kawa"lek
p"laszczyzny (bo sfera, na kt"or"a przeprowadza nas Gauss,
w okolicach punktu $(1,0,0)$ wygl"ada jak p"laszczyzna $z=1$).
Odwzorowanie to ma macierz $\left[ \begin{array}{cc} -2a & -2b \\ -2b & -2c \end{array} \right]$.
Czynnik, o jaki powi"eksza si"e pole powierzchni
przy odwzorowaniu Gaussa, jest r"owny wyznacznikowi tego przekszta"lcenia,
czyli $4(ac-b^2)$.

Je"sli por"ownasz ten wynik z wzorami na krzywizny g"l"owne, kt"ore
dostali"smy w rozdziale \ref{zewnetrzna}, oka"re si"e, "re jest on
r"owny iloczynowi krzywizn g"l"ownych.

Zatem {\bf krzywizna wewn"etrzna dowolnej powierzchni zanurzonej w przestrzeni tr"ojwymiarowej
jest r"owna iloczynowi krzywizn g"l"ownych}. Ten wynik Gauss nazwa"l {\em theorema egregum},
tak bardzo ucieszy"l go ten wynik. Czy ciebie te"r ucieszy"l?

\section{Parkiety.}
Rysujemy parkiet. Najpierw musimy si"e zdecydowa"c z czego
b"edziemy budowa"c nasz parkiet. Mo"remy si"e zdecydowa"c na
dowolny n-k"at: czworok"at, pi"eciok"at,...
Na pocz"atek um"owmy si"e,
"re b"edzie to tr"ojk"at.
Nast"epnie musimy si"e zdecydowa"c ile klepek parkietu ma si"e
styka"c w ka"rdym wierzcho"lku.
Zobaczmy kilka przyk"lad"ow...

...tu w ka"rdym wierzcho"lku styka si"e pi"e"c klepek. No c"o"r, nasz
parkiet jest do"s"c ko"slawy, ale przecie"r nikt nie obiecywa"l, "re parkiet
b"edzie zrobiony z wielok"at"ow foremnych.

...a tu w ka"rdym wierzcho"lku styka si"e siedem klepek.

Narysuj kilka tr"ojk"atnych parkiet"ow w kt"orych w ka"rdym
wierzcho"lku stykaja si"e: trzy, cztery, pi"e"c, sze"s"c i siedem klepek.

Jak to si"e dzieje, "re dla pewnych liczb klepek w wierzcho"lku (6 i
7) mo"rna zbudowa"c parkiet z dowolnie wielk"a liczb"a klepek, a dla
innych--nie? Co to ma wsp"olnego z charakterystyk"a Eulera?

Popatrz na twoje parkiety jak na sp"laszczone i zdeformowane
modele wielo"scian"ow foremnych. W ten spos"ob istota "ryj"aca w
dwuwymiarowym "swiecie mog"laby dok"ladnie policzy"c ilo"s"c
wierzcho"lk"ow, "scian i kraw"edzi dla wszystkich wielo"scian"ow
foremnych. Zr"ob to i ty!
Taka interpretacja parkietu za"lamuje si"e np. dla tr"ojk"atnego
parkietu o sze"sciu (lub wi"ecej) klepkach w wierzcho"lku. Jak wi"ec nale"ry
zinterpretowa"c taki parkiet?

Parkiety mog"lyby pom"oc istocie dwuwymiarowej w opisaniu
wielo"scian"ow foremnych. W jaki spos"ob my, istoty tr"ojwymiarowe
mo"remy zabra"c si"e za czterowymiarowe wielo"sciany foremne? (zerknij do ksi"a"rki Coxetera)

Obejrzyj (ale gdzie?) pi"ekne grafiki Eschera zatytu"lowane przez niego
{\em Kr"egi graniczne}. Co maj"a one wsp"olnego z naszymi parkietami?

\subsection{W"edr"owki po parkietach.}
W pewnym szalonym mie"scie plan miasta wygl"ada dok"ladnie tak,
jak jeden z naszych parkiet"ow. Wybieramy si"e na spacer.
Na ka"rdym skrzy"rowaniu zapisujemy, w kt"or"a drog"e skr"ecili"smy
W ten spos"ob otrzymujemy pewien napis okre"slaj"acy w
pe"lni nasz spacer, np: 224356 oznacza, "re dwa razy skr"ecili"smy w
ulic"e drug"a z lewej, potem w czwart"a... Mo"re wymy"slisz jaki"s
lepszy spos"ob?

W"edruj"ac po mie"scie nie zawsze robimy to w spos"ob optymalny,
tj. nie chodzimy najkr"otsz"a mo"rliw"a drog"a. Opisz, jak z tekstu
opisuj"acego pewien spacer otrzyma"c tekst opisuj"acy najkr"otszy
mo"rliwy spacer mi"edzy tymi samymi dwoma punktami. Jak nie robi"ac
rysunku udowodnisz, "re jest to naprawd"a najkr"otszy taki spacer?

Dwie osoby wybra"ly si"e na przechadzk"e. Znasz napisy
okre"slaj"ace ich spacery. Jak sprawdzisz, czy dosz"ly one do tego
samego miejsca?

Pewna osoba wybra"la si"e na w"edr"owk"e po trasie, kt"ora okaza"la
si"e by"c "laman"a zamkni"et"a. Jak zmierzy"c pole zawarte wewn"atrz tej
"lamanej, je"sli znasz napis okre"slaj"acy ten spacer? Przez pole
powierzchni umawiamy si"e rozumie"c ilo"s"c klepek parkietowych.

\subsection{Geometria na parkiecie.}
Wybierz sobie jaki"s parkiet. Policz ilo"s"c skrzy"rowa"n, do kt"orych
mo"rna doj"s"c w N krokach z ustalonego miejsca i zr"ob wykres.
Czym r"o"rni"a si"e przypadki parkietu tr"ojk"atnego w kt"orym w
wierzcho"lkach styka si"e pi"e"c, sze"s"c i siedem klepek?

Jak zale"ry pole ko"la od promienia na sferze, p"laszczy"znie i na
p"laszczy"znie hiperbolicznej?

\section{A co jest w wi"ekszych wymiarach?}
Nasze dotychczasowe rozwa"rania dotyzy"ly r\.o"rnych
dwuwymiarowych powierzchni. Ale oczywi"scie nie ma powodu, by nie
spr"obowa"c tego samego z zakrzywionymi powierzchniami tr"ojwymiarowymi,
czterowymiarowymi,... Czy co"s si"e zmieni, czy wszystko, co
powiedzieli"smy pozostaje prawd"a?

Poka"r, "re suma k"at"ow bry"lowych czworo"scianu nie jest zawsze
taka sama. A zatem nie mo"rna ju"r zdefiniowa"c jakiego"s
tr"ojwymiarowego odpowiednika defektu.

Udowodnij, "re charakterystyka Eulera dowolnej rozmaito"sci tr"ojwymiarowej
bez brzegu r"owna jest 0. Jak widzisz, nie jest "latwo uog"olni"c
wz"or Gaussa--Bonneta na wy"rsze wymiary (ale jest to mo"rliwe).

Osobi"scie podejrzewam, "re zasadnicza r"o"rnica pomi"edzy dwoma oraz
trzema lub wi"ecej wymiarami le"ry w tym, "re na dwuwymiarowej
p"laszczyznie mo"rna dobrze zdefiniowa"c k"at obrotu, a w wy"rszych
wymiarach nie (m"owi"ac ezoterycznie: grupa obrot"ow p"laszczyzny
$SO(2)$ jest przemnienna, a $SO(3)$ ju"r nie).

Je"sli chcesz wiedzie"c, jak mimo wszystko da"c sobie rad"e w
przestrzeni tr"ojwymiarowej (marzycielu, obud"z si"e! przyszl{}o ci
"ry"c w"la"snie w takiej przestrzeni!) si"egnij po ksi"a"rk"e
Feynmanna.

\section{Zrozumie"c Einsteina.}
Nie mam si"ly. Poczytaj Feynmanna. Albo zapytaj mnie.
\section{Odpowiedzi na niekt"ore zadania.}
\rozwiazanie{\ref{parabola}}
{\em Rozwi"azanie Tomasza Elsnera.}

Parabola o r"ownaniu $y=a x^2$ to zbi"or punkt"ow, kt"orych odleg"lo"s"c
od prostej $m$ o r"ownaniu $y=-l$ jest r"owna odleg"lo"sci od punktu $P=(0,l)$, gdzie
$a=\frac{1}{4 l}$ (sprawd"z!).

Je"sli $Q$ jest punktem naszej paraboli, a $R$ jest rzutem punktu $P$ na prost"a $m$,
to dwusieczna k"ata $\angle PQR$ jest styczn"a do paraboli wyznaczon"a
w punkcie $Q$. Gdyby bowiem tak nie by"lo, nasza dwusieczna przecina"laby
parabol"e w jakim"s innym punkcie $Q'$, kt"orego rzut na prost"a $m$.
Tr"ojk"aty $PQQ'$ i $RQQ'$
s"a przystaj"ace, poniewa"r maj"a wsp"olny bok $QQ'$, $PQ=QR$ oraz
k"aty $\angle PQQ'=\angle RQQ'$, zatem r"ownie"r $PQ'=Q'R$. Ale poniewa"r
punkt $Q'$ le"ry na paraboli, zatem $PQ'=Q'R'$. \L{}"acz"ac te wyniki
dostajemy $Q'R'=Q'R$, co  nie jest mo"rliwe.

Wybieramy punkt $Q$ tak, by by"l bardzo blisko punktu $O=(0,0)$, w kt"orym
to punkcie chcemy zmierzy"c krzywizn"e.
Narysujmy normalne do naszej paraboli w punktach $O$ i $Q$.
Oznaczmy ich punkt przeci"ecia przez $S$. Wida"c, "re proste $SQ$ i
$PR$ s"a prostopad"le do stycznej do paraboli wyznaczonej w punkcie
$Q$ czyli dwusiecznej k"ata $PQR$, a wi"ec te proste s"a r"ownoleg"le.
Punkty $SQRP$ tworz"a r"ownoleg"lobok, zatem $SP=QR$.

Promie"n krzywizny paraboli w punkcie $O$ zdefiniowali"smy jako odleg"lo"s"c
$SO=l+QR$. Je"sli punkt $Q$ wybrali"smy bardzo blisko $O$, to
$QR\approx PO=l$, zatem $SO=2l=\frac{1}{2a}$.

Promie"n krzywizny elipy na ko"ncu d"lu"rszej p"o"losi wynosi
$\frac{(a-b)(a+b)}{a}$, a na ko"ncu kr"otszej $\frac{a^2}{\sqrt{a^2-b^2}}$

\rozwiazanie{\ref{calkowita}}
Narysujmy w ka"rdym punkcie krzywej wektor styczny. Po przej"sciu "luku
okr"egu o d"lugo"sci $l$ i promieni krzywizny $r$, wektor styczny
obraca si"e o k"at r"owny $\frac{l}{r}$. Zatem po przej"sciu ca"lej
krzywej wektor styczny obr"oci si"e o k"at
$\frac{l_1}{r_1}+\frac{l_2}{r_2}+...$. Ale ca"lkowit"a krzywizn"e
krzywej zdefiniowali"smy jako t"e sum"e. Zatem ca"lkowita krzywizna
krzywej jest r"owna $2\pi\times$ (ile razy obraca si"e wektor styczny przy jednokrotnym
obej"sciu krzywej).

W szczeg"olno"sci ca"lkowita krzywizna okr"egu jest r"owna
$\pm2\pi$ (to, czy nale"ry wybra"c plus czy minus, zale"ry od tego,
w kt"or"a stron"e obchodzimy nasz okr"ag),
ca"lkowita krzywizna "osemki jest r"owna $0$...

\rozwiazanie{\ref{twkosinusow}}
$\cos c=\cos a \cos b+\sin a \sin b \cos \gamma $

\rozwiazanie{\ref{trygonometria}}
{\bf Rozwi"azanie 1.} Wybieramy uk"lad wsp"o"lrz"ednych tak samo,
jak w dowodzie twierdzenia Pitagorasa i teraz te"r $A=(\cos b,\sin b,0)$,
$B=(\cos a,0,\sin a)$, $C=(1,0,0)$. Szukamy teraz wielkiego ko"la przechodz"acego
przez $A$ i $B$. Nie jest to trudne: przecie"r wektor prostopad"ly do
wektor"ow $\stackrel{\longrightarrow}{OA}$ i $\stackrel{\longrightarrow}{OB}$
mo"rna otrzyma"c wyznaczaj"ac ich iloczyn wektorowy, kt"ory w tym wypadku
jest r"owny $(\sin a \sin b,-\cos a \sin b,-\sin a \cos b)$. Ten wektor ma
d"lugo"s"c r"own"a iloczynowi d"lugo"sci mno"ronych wektor"ow oraz
sinusa k"ata zawartego mi"edzy nimi, zatem nasz wektor ma d"lugo"s"c $\sin c$.
Zatem wektorem o d"lugo"sci $1$ prostopad"lym do $\stackrel{\longrightarrow}{OA}$ i $\stackrel{\longrightarrow}{OB}$
jest $(\frac{\sin a \sin b}{\sin c},-\frac{\cos a \sin b}{\sin c},-\frac{\sin a \cos b}{\sin c})        $ i
s"a to jednocze"snie wsp"o"lczynniki wielkiego ko"la wyznaczonego
przez bok $c$. Wielkie ko"lo wyznaczone przez bok $b$ ma oczywi"scie
wsp"o"lrz"edne $(0,1,0)$ lub $(0,-1,0)$. Wybiermy t"e drug"a
mo"rliwo"s"c--gdyby"smy wybrali pierwsz"a nie sta"loby si"e nic
strasznego, ale zamiast kosinusa k"ata $\alpha$ otrzymaliby"smy
kosinus k"ata dope"lniaj"acego do $\alpha$, czyli $\pi-\alpha$.
Tak wi"ec
$\cos \alpha=(0,-1,0)\cdot(\frac{\sin a \sin b}{\sin c},-\frac{\cos a \sin b}{\sin c},-\frac{\sin a \cos b}{\sin c})=
 \frac{\cos a \sin b}{\sin c}$. Z twierdzenia Pitagorasa $\cos a=\frac{\cos c}{\cos b}$,
zatem $\cos \alpha=\frac{\tan b}{\tan c}$.

{\bf Rozwi"azanie 2.} Ze wzor"ow al-Battaniego mamy
$$\cos \alpha=\frac{\cos a-\cos b \cos c}{\sin b \sin c}$$
...podnosimy obie strony do kwadratu...
$$\cos^2 \alpha=\frac{(\cos a-\cos b \cos c)^2}{\sin^2 b \sin^2 c}=$$
$$=\frac{(\cos a-\cos b \cos c)^2}{(1-\cos^2 b) \sin^2 c}=
\frac{(\cos^2 a-\cos a \cos b \cos c)^2}{(\cos^2 a-\cos^2 a \cos^2 b) \sin^2 c}$$
...teraz zauwa"ramy, "re z twierdzenia Pitagorasa $\cos a \cos b=\cos c$...
$$\cos^2 \alpha=\frac{(\cos^2 a-\cos^2 c)^2}{(\cos^2 a-\cos^2 c) \sin^2 c}=
\frac{\cos^2 a-\cos^2 c}{\sin^2 c}$$
$$\sin^2 \alpha=1-\cos^2 \alpha=\frac{\sin^2 a}{\sin^2 c}$$
A zatem mamy
$$\sin \alpha=\frac{\sin a}{\sin c}$$
(troch"e oszukali"smy, co ze znakiem?)

Podobnie znajdujemy drug"a formu"l"e:
$$\cos \alpha=\frac{\cos a-\cos b \cos c}{\sin b \sin c}=
\frac{\cos a \cos b-\cos^2 b \cos c}{\sin b \cos b \sin c}=$$
$$=\frac{\cos c-\cos^2 b \cos c}{\sin b \cos b \sin c}=
\frac{\sin b \cos c}{\cos b \sin c}=\frac{\tan b}{\tan c}$$

\rozwiazanie{\ref{obwod}}
Obw"od jest r"owny $2 \pi \sin r$, a pole r"owne $2 \pi (1-\cos r)$.

\rozwiazanie{\ref{promien}}
Odpowied"z 1:
Je"sli chodzimy tak, "reby jedna noga st"apa"la po
okr"egu o promieniu $a$, a druga po okr"egu o promieniu $a+\Delta a$,
to jedna noga przechodzi drog"e $2 \pi \sin a$, a druga $2 \pi \sin(a+\Delta a)=
2 \pi(\sin a \cos \Delta a+\cos a \sin \Delta a)\approx 2 \pi (\sin a+\Delta a \cos a)$.
St"ad znajdujemy promie"n krzywizny r"owny $\tan a$.

Odpowied"z 2: Rozwa"rmy sto"rek styczny do sfery wzd"lu"r naszego
okr"egu. Mieszka"ncy sfery i sto"rka mierz"ac promie"n krzywizny
naszego okr"egu otrzymaj"a ten sam wynik. A wynik, jaki otrzymaj"a
mieszka"ncy sto"rka "latwo ot\-rzy\-ma"c, poniewa"r (po drobnej ingerencji
no"ryczek) sto"rek daje si"e rozprostowa"c bez "radnych deformacji
do kawa"lka p"laszczyzny.

\rozwiazanie{\ref{zewnetrzna}}
Krzywizny g\l"owne naszej powierzchni s"a r"owne
$k_{min}=a+c-\sqrt{(a-c)^2+4 b^2}$
$k_{max}=a+c+\sqrt{(a-c)^2+4 b^2}$.
$k_{min}+k_{max}=2(a+c)$, $k_{min} k_{max}=4(ac-b^2)$.

\rozwiazanie{\ref{spacery}}
{\large \bf Spacer po okr"egu.}
Rozwa"rmy okr"ag na sferze (nie koniecznie ko"lo wielkie)--to po nim
b"edziemy w"edrowa"c. "Srodek sfery tradycyjnie oznaczamy przez
$O$, "srodek naszego ok\-r"e\-gu oznaczmy przez $R$. Dla ka"rdego punktu $P$
le"r"acego na okr"egu k"at $\angle POR$ jest taki sam, oznaczmy
wi"ec go przez $\alpha$.

Promie"n naszego okr"egu jest r"owny oczywi"scie $\sin \alpha$.
Tym nie mniej mie\-szka\-"nco\-wi sfery wydaje si"e, "re ma on promie"n
$\tan \alpha$. Na p"laszczy"znie trzeba przej"s"c d"lugo"s"c
$2 \pi \tan \alpha$ wzd"lu"r "luku okr"egu o tym promieniu, "reby
wr"oci"c do punktu wyj"scia. Tymczasem id"ac tyle samo wzd"lu"r
naszego okr"egu na sferze, wracamy do punktu wyj"scia, mijamy go i
idziemy jeszcze troch"e dalej. Faktyczny obw"od okr"egu jest r"owny
$2 \pi \sin \alpha$, zatem po mini"eciu punktu pocz"atkowego przejdziemy
jeszcze $2 \pi \tan \alpha-2 \pi \sin \alpha=2 \pi \tan \alpha (1-\cos \alpha)$.

Cho"c mogliby"smy dok"ladnie policzy"c odleg"lo"s"c punktu pocz"atkowego
i ko"ncowego naszej podr"o"ry, pos"lu"rymy si"e prostym przybli"rniem.
Mianowicie w przybli"reniu ta odleg"lo"s"c jest r"owna d"lugo"sci
drogi, jak"a przechodzimy po mini"eciu punktu pocz"atkowego, czyli
$2 \pi \tan \alpha (1-\cos \alpha)\approx \pi \alpha^3 $.

A teraz zastan"owmy si"e, jaki k"at tworz"a pocz"atkowy i ko"ncowy
kierunek naszego nosa. I tu pojawia si"e problem: mianowicie chcemy mierzy"c
k"at pomi"edzy dwoma wektorami (nosy), kt"orych punkty zaczepienia
s"a gdzie indziej. Fachowiec od geometrii r"o"rniczkowej stanowczo
zaprotestowa"lby przeciwko takim pr"obom.

Rozwa"rmy wi"ec nieco zmodyfikowany spacer na p"laszczy\,znie.
B"edziemy chodzi"c (tak jak poprzednio) po okr"egu o promieniu $R=\tan \alpha$,
ale tym razem nie obejdziemy ca"lego obwodu, ale tylko d"lugo"s"c
$2 \pi \sin \alpha$. Zatem k"at mi"edzy pocz"atkowym i
ko"ncowym kierunkiem nosa jest r"owny storunkowi przebytej drogi i
promienia ko"la, czyli $\frac{2 \pi \sin \alpha}{\tan \alpha}=2 \pi \cos \alpha$.
Je"sli  wi"ec teraz obr"ocimy si"e jeszcze o k"at $2 \pi (1-\cos \alpha)$,
to, tak jak chcieli"smy, nos wr"oci do pocz"atkowego po"lo"renia.

Je"sli ten sam spacer wykonamy na sferze, to zgodnie z tym, co ju"r powiedzieli"smy,
wr"ocimy dok"ladnie do punktu startu; ale obr"ot o $2 \pi (1-\cos \alpha)$
sprawia, "re w"la"snie o tyle r"o"rni si"e ko"ncowe i pocz"atkowe
po"lo"renie nosa. W przybli"reniu jest to r"owne $\pi \alpha^2$.

{\em Wynik ten nie powinien nas dziwi"c. Przecie"r k"at, o jaki zostajemy
o\-br"o\-ce\-ni jest dok"ladnie r"owny defektowi wielok"ata po kt"orym
idziemy. Na sferze o promieniu $1$ defekt jest dok"ladnie r"owny polu,
a pole czaszy, kt"or"a obchodzimy jest przecie"r r"owne
$2 \pi (1-\cos \alpha)$.}

Wypada teraz wyci"agn"a"c wnioski. W"edruj"ac po sferze napotykamy
dwa zjawiska wywo"lane przez krzywizn"e: nie wracamy dok"ladnie do
punktu startu (cho"c 'powinni"smy') i zostajemy troch"e obr"oceni.
Wielko"s"c ka"rdego z tych efekt"ow zale"ry--m"owi"ac troch"e
nieprecyzyjnie--od wielko"sci obszaru, po kt"o\-rym si"e poruszamy.
W przypadku chodzenia po okr"egu mo"remy t"e wielko"s"c scharakteryzowa"c
np. przez $\alpha$. Pierwszy efekt jest rz"edu $\alpha^3$, czyli jest
dla ma"lych $\alpha$ pomijalnie ma"ly w por"ownaniu z wielko"sci"a
drugiego efektu, kt"ory jest r"owny polu obchozonego obszary (w przybli"reniu $\alpha^2$).

{\large \bf Spacer po prostok"acie.}
Idziesz prosto przed siebie (przebyt"a drog"e oznaczymy przez $a$).
Nast"epnie idziesz 'rakiem' w prawo i przebywasz drog"e $b$. Teraz
idziesz do ty"lu, przebywasz drog"e $a$. I znowu idziesz rakiem, ale w
lewo i przebywasz drog"e $b$. Gdzie jeste"s?

Na p"laszczy"znie odpowiedz jest prosta: jestem dok"ladnie tam, gdzie
by"lem na samym pocz"atku. A na sferze?

Wprowad"zmy uk"lad wsp"o"lrz"ednych sztywno zwi"azany z w"edrowcem.
W"edr"owka w"edrowca po sferze wygl"ada w takim uk"ladzie troch"e
inaczej: w"edrowiec spoczywa, a pod jego stopami obraca si"e sfera.
Ka"rdemu takiemu obrotowi sfery odpowiada jaka"s macierz.
Macierze odpowiadaj"ace kolejnym krokom s"a r"owne
$\left[\begin{array}{ccc} \cos a & -\sin a & 0 \\ \sin a & \cos a & 0 \\
0 & 0 & 1 \end{array} \right]$, $\left[\begin{array}{ccc} 1 & 0 & 0 \\
0 & \cos b & -\sin b \\ 0 & \sin b & \cos b \end{array} \right]$,
$\left[\begin{array}{ccc} \cos a & \sin a & 0 \\ -\sin a & \cos a & 0 \\
0 & 0 & 1 \end{array} \right]$, $\left[\begin{array}{ccc} 1 & 0 & 0 \\
0 & \cos b & \sin b \\ 0 & -\sin b & \cos b \end{array} \right]$. \.Zeby
zobaczy"c macierz odpowiadaj"ac"a ca"lemu spacerowi, trzeba je przemno"ry"c
w tej w"la"snie kolejno"sci. Najpierw policzymy iloczyn pierwszych dwu:
$$ \left[\begin{array}{ccc} \cos a & -\sin a & 0 \\ \sin a & \cos a & 0 \\
0 & 0 & 1 \end{array} \right] \times \left[\begin{array}{ccc} 1 & 0 & 0 \\
0 & \cos b & -\sin b \\ 0 & \sin b & \cos b \end{array} \right] =
\left[ \begin{array}{ccc} \cos a & -\sin a \cos b & \sin a \sin b \\
\sin a & \cos a \cos b & -\cos a \sin b \\ 0 & \sin b & \cos b \end{array} \right] $$
Poniewa"r macierz trzecia i czwarta r"o"rni"a sie od pierwszej i drugiej
tylko znakami przy sinusach, iloczyn trzeciej i czwartej macierzy "latwo
uzyska"c zmieniaj"ac znaki przy sinusach w obliczonym ju"r iloczynie
pierwszej i drugiej macierzy. Tak wi"e iloczyn wszystkich czterech macierzy
jest r"owny
$$\left[ \begin{array}{ccc} \cos a & -\sin a \cos b & \sin a \sin b \\
\sin a & \cos a \cos b & -\cos a \sin b \\ 0 & \sin b & \cos b \end{array} \right]
\times \left[ \begin{array}{ccc} \cos a & \sin a \cos b & \sin a \sin b \\
-\sin a & \cos a \cos b & \cos a \sin b \\ 0 & -\sin b & \cos b \end{array} \right]= $$
$$\left[ \begin{array}{ccc} 1-\sin^2 a (1-\cos b)  
& \sin a \cos a \cos b (1-\cos b) & \sin a \sin b [1-(1-\cos a)(1-\cos b)] 
\\ \sin a \cos b (1-\cos a) & 1+(1-\cos a)(1-\cos b)(\cos a \cos b-1) &
\sin b [\sin^2 a+\cos a \cos b (cos a-1)]  \\
-\sin a \cos b & \sin b \cos b (\cos a-1) & 1-\sin^2 b (1-\cos a) \end{array} \right] $$
$$\approx \left[ \begin{array}{ccc} 1-\frac{a^2 b^2}{2} & \frac{a b^2}{2} & a b \\
\frac{a^3}{2} & 1-\frac{a^2 b^2 (a^2+b^2)}{8} & \frac{a^2 b}{2} \\
-a b & -\frac{a^2 b}{2} & 1-\frac{a^2 b^2}{2} \end{array} \right]
$$
Wnioski, jakie mo"rna wyci"agn"a"c z tego wyniku s"a takie same,
jak dla w"edr"owki po okr"egu.

\end{document}